\documentclass[a4paper, DIV=12, 11pt]{scrartcl}
\usepackage[latin1]{inputenc}
\usepackage{amsmath}
\usepackage{amsfonts}
\usepackage{amssymb}
\usepackage{amsthm}
\usepackage{graphicx}
\usepackage{thmtools}
\usepackage{bbm}
\usepackage{hyperref}
\usepackage[bottom, marginal]{footmisc}
\usepackage{todonotes}

\usepackage{tikz}
\usetikzlibrary{matrix}

\declaretheoremstyle[headfont=\normalfont]{normalhead}
\newtheorem{lemma}{Lemma}[section]
\newtheorem{theorem}[lemma]{Theorem}
\newtheorem{proposition}[lemma]{Proposition}
\newtheorem{corollary}[lemma]{Corollary}

\newcommand{\R}{\mathbb{R}}
\newcommand{\C}{\mathbb{C}}

\DeclareMathOperator{\sign}{sign}

\DeclareMathOperator{\VConv}{VConv}
\DeclareMathOperator{\Conv}{Conv}
\DeclareMathOperator{\vol}{vol}

\DeclareMathOperator{\supp}{supp}

\DeclareMathOperator{\D}{\bar{D}}

\DeclareMathOperator{\GW}{GW}
\DeclareMathOperator{\SO}{\mathrm{SO}}

\DeclareMathOperator{\nc}{\mathrm{nc}}

\setkomafont{title}{\normalfont\Large}

\author{Jonas Knoerr}
\title{Singular valuations and the Hadwiger theorem on convex functions}
\date{}

\newcommand{\Addresses}{{
		\bigskip
		\footnotesize
		
		Jonas Knoerr, \textsc{Institute of Discrete Mathematics and Geometry, TU Wien, Wiedner Hauptstrasse 8-10, 1040 Wien, Austria}\par\nopagebreak
		\textit{E-mail address}: \texttt{jonas.knoerr@tuwien.ac.at}
		
		\medskip
	}}
	
\makeatletter
\def\blfootnote{\xdef\@thefnmark{}\@footnotetext}
\makeatother

\makeindex
\begin{document}
\maketitle
\begin{abstract}
	We give a characterization of smooth, rotation and dually epi-translation invariant valuations and use this result to obtain a new proof of the Hadwiger theorem on convex functions, originally established by Colesanti, Ludwig and Mussnig. We also give a description of the construction of the functional intrinsic volumes using integration over the differential cycle and provide a new representation of these functionals as principal value integrals with respect to the Hessian measures.
\end{abstract}	

\section{Introduction}

Given a family $X$ of (extended) real-valued functions, a functional $\mu:X\rightarrow\R$ is called a valuation if
\begin{align*}
	\mu(f\vee h)+\mu(f\wedge h)=\mu(f)+\mu(h)
\end{align*}
for all $f,h\in X$ such that the pointwise maximum $f\vee h$ and minimum $f\wedge h$ belong to $X$. If $X$ denotes the class of indicator functions associated to convex bodies of $\R^n$, that is, the family of non-empty convex and compact subsets of $\R^n$, this recovers the classical notion of valuations on convex bodies.\\

In recent years, many classical results about valuations on convex bodies were generalized to valuations on different classical function spaces, for example $L_p$ \cite{TsangValuations$Lp$spaces2010} and Sobolev spaces \cite{LudwigFisherinformationmatrix2011,LudwigValuationsSobolevspaces2012,LudwigCovariancematricesvaluations2013}, Lipschitz functions \cite{ColesantiEtAlclassinvariantvaluations2020,ColesantiEtAlContinuousvaluationsspace2021}, definable functions \cite{BaryshnikovEtAlHadwigersTheoremdefinable2013}, or general Banach lattices \cite{TradaceteVillanuevaValuationsBanachlattices2020}.

In a series of articles, Colesanti, Ludwig and Mussnig proved a version of the Hadwiger theorem for valuations on convex functions and provided a number of descriptions of the relevant valuations \cite{ColesantiEtAlHadwigertheoremconvex,ColesantiEtAlHadwigertheoremconvex2020,ColesantiEtAlHadwigertheoremconvex2022,ColesantiEtAlHadwigertheoremconvex2023}. Let us introduce the framework of this article and the version of their result we will be interested in. Let $\Conv(\R^n,\R)$ denote the space of finite-valued convex functions on $\R^n$, which is a metrizable space when equipped with the topology induced by epi-convergence (which coincides with pointwise convergence in this setting, compare Section \ref{section:convexFunctions}). We will be interested in \emph{dually epi-translation invariant} valuations, that is, valuations $\mu:\Conv(\R^n,\R)\rightarrow\R$ that satisfy 
\begin{align*}
	\mu(f+\lambda+c)=\mu(f)\quad\text{for all } f\in\Conv(\R^n,\R), \lambda:\R^n\rightarrow\R \text{ linear},c\in \R. 
\end{align*}
This notion is intimately tied to translation invariance and thus there is a strong connection with the theory of translation invariant valuations on convex bodies, see \cite{ColesantiEtAlhomogeneousdecompositiontheorem2020,KnoerrSmoothvaluationsconvex2024}. A weaker notion is considered in \cite{AleskerValuationsconvexfunctions2019}.\\
Let $\VConv(\R^n)$ denote the space of all continuous, dually epi-translation invariant valuations on $\Conv(\R^n,\R)$ equipped with the topology of uniform convergence on compact subsets. As shown in \cite{ColesantiEtAlhomogeneousdecompositiontheorem2020}, there exists a homogeneous decomposition of this space: If we denote by $\VConv_i(\R^n)\subset\VConv(\R^n)$ the subspace of all $i$-homogeneous valuations, that is, all valuations $\mu$ that satisfy $\mu(tf)=t^i\mu(f)$ for $t\ge 0$, then
\begin{align*}
	\VConv(\R^n)=\bigoplus_{i=0}^n\VConv_i(\R^n).
\end{align*}
Obviously, we have a natural action of the special orthogonal group $\SO(n)$ on $\VConv(\R^n)$ given by $[g\cdot\mu](f):=\mu(f\circ g)$ and this action respects the homogeneous decomposition. The Hadwiger theorem on convex functions provides a characterization of the space $\VConv_i(\R^n)^{\SO(n)}$ of all $\SO(n)$-invariant valuations in $\VConv_i(\R^n)$. To state it, we first need to discuss the relevant class of invariant valuations introduced by Colesanti, Ludwig and Mussnig in \cite{ColesantiEtAlHadwigertheoremconvex2020}. Let $C_b((0,\infty))$  denote the space of all continuous functions on $(0,\infty)$ with support bounded from above and consider the subspace $D^n_i$ for $1\le i\le n-1$ given by
\begin{align*}
	D^n_i:=\left\{\zeta\in C_b((0,\infty)):\lim\limits_{t\rightarrow0}t^{n-i}\zeta(t)=0, \lim\limits_{t\rightarrow0}\int_t^\infty \zeta(s)s^{n-i-1}ds\text{ exists and is finite}\right\}.
\end{align*}
\begin{theorem}[\cite{ColesantiEtAlHadwigertheoremconvex2020} Theorem 1.4]
	\label{theorem:existenceFunctionalIntrinsicVolumes}
	Let $1\le i\le n-1$. For every $\zeta\in D^n_i$ there exists a unique valuation $V^*_{i,\zeta}\in \VConv_i(\R^n)^{\SO(n)}$ such that 
	\begin{align*}
		V^*_{i,\zeta}(f)=\int_{\R^n}\zeta(|x|)[D^2f(x)]_i dx
	\end{align*}
	for all $f\in\Conv(\R^n,\R)\cap C^2(\R^n)$.
\end{theorem}
Here, $[D^2f(x)]_i$ denotes the $i$th elementary symmetric polynomial in the eigenvalues of the Hessian $D^2f(x)$ of $f$ in $x\in\R^n$. \\
Using these valuations, Colesanti, Ludwig and Mussnig established the following complete characterization of the spaces $\VConv_i(\R^n)^{\SO(n)}$:
\begin{theorem}[\cite{ColesantiEtAlHadwigertheoremconvex2020}]
	\label{maintheorem:continuousHadwiger}
	Let $1\le i\le n-1$. For every $\mu\in\VConv_i(\R^n)^{\SO(n)}$ there exists a unique $\zeta\in D^n_i$ such that $\mu=V^*_{i,\zeta}$.
\end{theorem}
Let us remark that we have excluded the cases $i=0$ and $i=n$ from the statement as the classes of functions are slightly different and the classification is much simpler. However, all of the relevant classes are related in a sensible manner, see \cite{ColesantiEtAlHadwigertheoremconvex2022}. Due to the parallel of this result to Hadwiger's classical classification of all rotation invariant, continuous and translation invariant valuations on convex bodies as linear combination of intrinsic volumes, the functionals $V^*_{i,\zeta}$, $\zeta\in D^n_i$ are called \emph{functional intrinsic volumes}. Let us also remark that Colesanti, Ludwig and Mussnig established an alternative representation of  the functional intrinsic volumes using averages with respect to restrictions to $i$-dimensional subspaces in \cite{ColesantiEtAlHadwigertheoremconvex2022}, which also provides a more direct representation of these valuations that does not involve singular integrals.\\

In \cite{ColesantiEtAlHadwigertheoremconvex2022}, Colesanti, Ludwig and Mussnig provided an alternative proof of Theorem \ref{maintheorem:continuousHadwiger} using the following classification of the dense subspace of so called \emph{smooth} rotation invariant valuations (see Section \ref{section:smoothValuations} for the definition), originally established by the author in his PhD thesis \cite{KnoerrSmoothvaluationsconvex2020a}:
\begin{theorem}[\cite{KnoerrSmoothvaluationsconvex2020a} Theorem 9.4.4]
	\label{theorem:smoothHadwiger}
	Let $1\le i\le n$. For every smooth valuation $\mu\in\VConv_i(\R^n)^{\SO(n)}$ there exists $\phi\in C^\infty_c([0,\infty))$ such that
	\begin{align*}
		\mu(f)=\int_{\R^n}\phi(|x|^2)[D^2f(x)]_idx
	\end{align*}
	for all $f\in\Conv(\R^n,\R)\cap C^2(\R^n)$.
\end{theorem}

The goal of this article is twofold: First, we provide a new proof of Theorem \ref{maintheorem:continuousHadwiger} using an approximation by smooth valuations similar to the proof given in \cite{ColesantiEtAlHadwigertheoremconvex2022}. For this, we also include a streamlined version of the proof of Theorem \ref{theorem:smoothHadwiger} from \cite{KnoerrSmoothvaluationsconvex2020a} in Section \ref{section:smoothHadwiger}. 
The second and main goal of this article is a reconceptualization of the original proof of Theorem \ref{theorem:existenceFunctionalIntrinsicVolumes} that can be generalized to the construction of "singular invariant valuations" for more general groups. This will also lead us to an interpretation of the functional intrinsic volumes as principal value integrals of the functions $\zeta(|\cdot|)$, $\zeta\in D^n_i$, as well as a complete answer to the question under which conditions on $\zeta\in D^n_i$ the functional intrinsic volume is given by the integral representation in Theorem \ref{theorem:existenceFunctionalIntrinsicVolumes} with respect to $f\in\Conv(\R ^n,\R)$.\\

Let us thus revisit the construction of the functional intrinsic volumes given by Colesanti, Ludwig and Mussnig in \cite{ColesantiEtAlHadwigertheoremconvex2020}. First, for continuous $\zeta\in C_c([0,\infty))$, $V^*_{i,\zeta}$ may be represented by
\begin{align*}
	V^*_{i,\zeta}(f)=\int_{\R^n}\zeta(|x|)d\Phi_i(f,x)\quad\text{for all }f\in\Conv(\R^n,\R),
\end{align*}
where $\Phi_i(f)$ denotes the so called $i$th Hessian measure of $f$. These measures were originally studied by Trudinger and Wang in \cite{TrudingerWangHessianmeasures.I1997,TrudingerWangHessianmeasures.II1999}. For their role in the construction of valuations see \cite{ColesantiHugSteinertypeformulae2000,ColesantiHugHessianmeasuresconvex2005,ColesantiEtAlHessianvaluations2020}. If $\zeta$ is not continuous in $0$, then $x\mapsto \zeta(|x|)$ is in general not integrable with respect to the measure $\Phi_i(f)$ and we are left with certain singular integrals, called \emph{singular Hessian valuations} in \cite{ColesantiEtAlHadwigertheoremconvex2020}. The proof of Theorem \ref{theorem:existenceFunctionalIntrinsicVolumes} given by Colesanti, Ludwig and Mussnig in \cite{ColesantiEtAlHadwigertheoremconvex2020} uses an elaborate regularization procedure to show that this expression extends  by continuity from smooth to arbitrary convex functions if $\zeta\in D^n_i$.\\
The main idea is to replace the integral above by a principal value type integral: If $f\in\Conv(\R^n,\R)$ is a $C^2$-function with uniform lower bound on its Hessian, then its Legendre transform $\mathcal{L}f$ has compact sublevel sets. If we assume in addition that $f(0)=0\le f(x)$ for all $x\in \R^n$, then the sublevel sets contain the origin in their interior and for $t>0$ the expression 
\begin{align*}
	\int_{\{x:\ t<\mathcal{L}f(d f(x))\}}\zeta(|x|)d\Phi_i(f,x)
\end{align*}
is therefore well defined. Here, we interpret $\mathcal{L}f$ as a function on the dual space $(\R^n)^*$ of $\R^n$ and use that the differential $df:\R^n\rightarrow(\R^n)^*$ defines a diffeomorphism between the two spaces. In fact, this expression is well defined for convex functions of lower regularity. For the construction of the functional intrinsic volumes, one then has to show that these quantities behave sufficiently well for $t\rightarrow0$, at least along certain sequences, and the limit along such sequences is taken as the value of $V^*_{i,\zeta}$. Thus one of the main problems in establishing the existence of the functional intrinsic volumes is deriving suitable estimates for integrals over such sublevel sets.\\

Note that for sufficiently regular functions $f\in\Conv(\R^n,\R)$, the functional intrinsic volumes can be interpreted as integrals over the graph of the differential $df$ of certain differential forms defined on the cotangent bundle. Our approach is based on the observation that all of the geometric steps in the construction of the functional intrinsic volumes can be reformulated in terms of certain natural operations on these forms. The differential forms we will use are precisely the forms one encounters in the classification of smooth rotation invariant valuations (compare Section \ref{section:smoothValuations} and Section \ref{section:InvDiffFormsAndHessiansMeasures}). \\

This allows us to reinterpret a version of the inequalities obtained in \cite{ColesantiEtAlHadwigertheoremconvex2020} for the sublevel sets as bounds on the operator norm of the map
\begin{align*}
	V^*_i:C_c([0,\infty))&\rightarrow\VConv_i(\R^n)^{\SO(n)}\\
	\zeta&\mapsto \left[f\mapsto \int_{\R^n}\zeta(|x|)d\Phi_i(f,x)\right]
\end{align*}
with respect to a certain norm on $C_c([0,\infty))$, and the unique continuous extension of this map to a suitable completion of $C_c([0,\infty))$ provides the functional intrinsic volumes. More precisely, for $R>0$ set 
\begin{align*}
	D^n_{i,R}:=\{\zeta\in D^n_i:\supp\zeta\subset(0,R]\}.
\end{align*}
We will see in Section \ref{section:Dnk} that $C_c([0,\infty))\cap  D^n_{i,R}$ is dense in $D^n_{i,R}$ with respect to the norm
\begin{align*}
	\|\zeta\|:=\sup_{t>0}\left|(n-i)\int_t^\infty \zeta(s)s^{n-i-1}ds\right|+\sup_{t>0}\left|t^{n-i}\zeta(t)+(n-i)\int_t^\infty \zeta(s)s^{n-i-1}ds\right|.
\end{align*}
Similarly, we consider the space $\VConv_{B_R(0)}(\R^n)$ of all valuations that are supported on the ball $B_R(0)$ with radius $R>0$ centered at the origin (see Section \ref{section:duallyTopologySupport} for the definition of this notion of support). As shown in \cite{Knoerrsupportduallyepi2021}, this is a Banach space, which allows us to give the following independent proof of the results in Theorem \ref{theorem:existenceFunctionalIntrinsicVolumes} and Theorem \ref{maintheorem:continuousHadwiger}.
\begin{theorem}
	\label{maintheorem:TopIsomorphism}
	Let $1\le i\le n-1$. For every $R>0$ the restriction of $V^*_i$ to $C_c([0,\infty))\cap D^n_{i,R}$ extends uniquely to a topological isomorphism
	\begin{align*}
		V^*_i:D^n_{i,R}\rightarrow\VConv_i(\R^n)^{\SO(n)}\cap \VConv_{B_R(0)}(\R^n).
	\end{align*}
\end{theorem} 
Let us add the following remark on the relation of map $V^*_i$ and the functional intrinsic volumes. By definition, $V^*_i(\zeta)$ coincides with the functional intrinsic volume $V^*_{i,\zeta}$ from Theorem \ref{theorem:existenceFunctionalIntrinsicVolumes} for $\zeta\in C_c([0,\infty))$. We will establish in Section \ref{Section:functionalIntrinsicVolumes} that this is in fact true for arbitrary $\zeta\in D^n_{i,R}$. Thus, Theorem \ref{maintheorem:TopIsomorphism} implies in particular that the functional intrinsic volume $V^*_{i,\zeta}$ depends continuously on $\zeta\in D^n_{i,R}$ with respect to the norm $\|\cdot\|$.\\

It will also become clear from the construction that the valuation $V^*_i(\zeta)$ does not depend on the choice of $R>0$ such that $\zeta\in D^n_{i,R}$, so we obtain a well defined map $V^*_i:D^n_{i}\rightarrow\VConv_i(\R^n)^{\SO(n)}$.\\
Note that by the main results of \cite{Knoerrsupportduallyepi2021} every valuation in $\VConv(\R^n)$ has compact support, so Theorem \ref{maintheorem:TopIsomorphism} shows that $V^*_i:D^n_i\rightarrow\VConv_i(\R^n)^{\SO(n)}$ is bijective, which implies Theorem \ref{maintheorem:continuousHadwiger}.\\

As an application of our estimates for $V^*_{i,\zeta}$, we establish the following new representation of the functional intrinsic volumes as principal value integrals with respect to the Hessian measures $\Phi_i$.
\begin{theorem}
	\label{maintheorem:PrincipalValue}
	For $\zeta\in D^n_i$ and $f\in\Conv(\R^n,\R)$,
	\begin{align*}
		V^*_{i,\zeta}(f)=\lim\limits_{\epsilon\rightarrow0}\int_{\R^n\setminus B_\epsilon(0)}\zeta(|x|)d\Phi_i(f,x).
	\end{align*}
	Moreover, the convergence is uniform on compact subsets in $\Conv(\R^n,\R)$. 
\end{theorem}
We refer to \cite{Knoerrsupportduallyepi2021} Proposition 2.4 for a descriptions of the compact subsets of $\Conv(\R^n,\R)$ with respect to the topology induced by epi-convergence. Let us note that Theorem \ref{maintheorem:PrincipalValue} has the following direct consequence: If the function $x\mapsto \zeta(x)$ is integrable with respect to $\Phi_i(f)$, then $V^*_{i,\zeta}(f)$ is given by the integral representation in Theorem \ref{theorem:existenceFunctionalIntrinsicVolumes} (see Corollary \ref{corollary:RepresentationFunctIntrinsicVolume} below). We provide a precise characterization when this is the case simultaneously for all $f\in\Conv(\R ^n,\R)$ in Corollary \ref{corollary:ZetaIntegrable}.

\subsection{Plan of this article}
In Section \ref{section:preliminaries} we recall some basic notions about smooth valuations on convex bodies and the representation of the intrinsic volumes in terms of integration of differential forms over the normal cycle. We also recall some basic facts about convex functions and discuss the classes $D^n_{i,R}$ and the norm defined above.\\
Section \ref{section:VConv} recalls some results about dually epi-translation invariant valuations and smooth valuations, which includes a key result needed in our proof of the Hadwiger theorem.\\
In Section \ref{section:PropDifferentialCycle} we discuss some continuity properties of the differential cycle and its relation to the normal cycle of sublevel sets for certain convex functions. The differential forms needed in the constructions and their relation to the Hessian measures are presented in Section \ref{section:InvDiffFormsAndHessiansMeasures}.\\
The construction of the functional intrinsic volumes is discussed in Section \ref{section:singularHessianValuations}, which also includes the proof of Theorem \ref{maintheorem:PrincipalValue}. Finally, we establish the classification of all smooth dually epi-translation invariant valuations in this class in Section \ref{section:Hadwiger} and show how one can extend this result to the continuous case, completing the proof of Theorem \ref{maintheorem:TopIsomorphism}.

\section{Preliminaries}
\label{section:preliminaries}
\subsection{Valuations on convex bodies and intrinsic volumes}
\label{section:IntegralGeo}
Let $\mathcal{K}(\R^n)$ denote the space of convex bodies in $\R^n$, that is, the space of all non-empty, compact and convex subsets of $\R^n$ equipped with the Hausdorff metric. A functional $\mu:\mathcal{K}(\R^n)\rightarrow\R$ is called a valuation if
\begin{align*}
	\mu(K\cup L)+\mu(K\cap L)=\mu(K)+\mu(L)
\end{align*}
for all $K,L\in\mathcal{K}(\R^n)$ such that $K\cup L$ is convex. We refer to the monograph by Schneider \cite{SchneiderConvexbodiesBrunn2014} for a background on convex bodies and valuations on convex bodies.\\

One way to obtain valuations on convex bodies is the following construction (for an extension of these notions to larger classes of sets see \cite{FuCurvaturemeasuressubanalytic1994} and the references therein):\\
For a convex body $K\in\mathcal{K}(\R^n)$, the set
\begin{align*}
	\nc(K):=\{(x,v)\in\R^n\times S^{n-1}:v\text{ outer normal to } K \text{ in }x\in\partial K\}
\end{align*}
is a compact Lipschitz submanifold of the sphere bundle $S\R^n:=\R^n\times S^{n-1}$ of dimension $n-1$ that carries a natural orientation. We may thus consider $\nc(K)$ as an integral $(n-1)$-current on $S\R^n$, called the \emph{normal cycle} of $K$, which we will denote by $\omega\mapsto \nc(K)[\omega]$, where $\omega\in\Omega^{n-1}(S\R^n)$ is a smooth $(n-1)$-form.

\medskip
Given $\omega\in \Omega^{n-1}(S\R^n)$, we may in particular associate to any $K\in\mathcal{K}(\R^n)$ the signed measure 
\begin{align*}
	B\mapsto \Phi_\omega(K,B):=\nc(K)[1_{\pi^{-1}(B)}\omega]\quad\text{for Borel sets }B\subset \R^n.
\end{align*}
Here $\pi:S\R^n\rightarrow\R^n$ denotes the projection onto the first factor and $1_A$ denotes the indicator function of a set $A$. The map $\Phi_\omega$ is called a smooth curvature measure.\\

Similarly, we obtain the real-valued functional $K\mapsto\nc(K)[\omega]$ on $\mathcal{K}(\R^n)$, which turns out to be a continuous valuation, see \cite{AleskerFuTheoryvaluationsmanifolds.2008} Theorem 2.1.13. The differential form inducing such a valuation is highly non-unique. The kernel of this assignment can be described using notions from contact geometry, which we briefly recall. $S\R^n$ is a contact manifold with $1$-form $\alpha$ given by $\alpha|_{(x,v)}(w)=\langle v,d\pi (w)\rangle$ for $w$ belonging to the tangent space $T_{(x,v)}S\R^n$ in $(x,v)\in S\R^n$, where $\pi:S\R ^n\rightarrow\R ^n$ denotes the natural projection. Moreover, one can show that for every $\omega\in\Omega^{n-1}(S\R^n)$ there exists a unique form $\alpha\wedge\xi$ such that $d(\omega+\alpha\wedge\xi)$ is \emph{vertical}, that is, a multiple of $\alpha$. The Rumin operator is then defined by $D\omega=d(\omega+\alpha\wedge \xi)$, compare \cite{RuminFormesdifferentiellessur1994}.

The following description of the kernel was obtained by Bernig and Br\"ocker in \cite{BernigBroeckerValuationsmanifoldsRumin2007}.
\begin{theorem}[\cite{BernigBroeckerValuationsmanifoldsRumin2007} Theorem 1]
	\label{theorem:KernelTheorem}
	Let $\omega\in\Omega^{n-1}(S\R^n)$ be given. 
	\begin{enumerate}
		\item $\nc(K)[\omega]=0$ for all $K\in\mathcal{K}(\R^n)$ if and only if
		\begin{enumerate}
			\item $D\omega=0$
			\item $\int_{\{x\}\times S^{n-1}}\omega=0$ for all $x\in\R^n$.
		\end{enumerate}
		\item $\Phi_\omega=0$ if and only if $\omega$ is contained in the ideal generated by $\alpha$ and $d\alpha$.
	\end{enumerate} 
\end{theorem}

We also recall the following notion: The unique vector field $T$ on $S\R^n$ with $i_T\alpha=1$, $i_Td\alpha=0$ is called the Reeb vector field. Note that the relation $\omega=\alpha\wedge i_T\omega$ holds for any vertical form.\\

We will call a smooth curvature measure $\Phi_\omega$ non-negative if $\Phi_\omega(K,\cdot)$ is a non-negative measure on $\R^n$ for all $K\in\mathcal{K}(\R^n)$. Similarly, we will call a valuation $\mu$ on $\mathcal{K}(\R^n)$
\begin{itemize}
	\item positive if $\mu(K)\ge 0$ for all $K\in\mathcal{K}(\R^n)$,
	\item monotone if $K\subset L$ implies $\mu(K)\le \mu(L)$.
\end{itemize}
Then the following holds:
\begin{lemma}[\cite{BernigFuHermitianintegralgeometry2011} Theorem 2.7]\label{lemma:monotoneValuationsCurvatureMeasures}
	$\omega\in\Omega^{n-1}(S\R^n)$ induces a monotone valuation if and only if $i_TD\omega$ induces a non-negative curvature measure.
\end{lemma}
Here, $i_X\omega$ denotes the interior product of a vector field $X$ and a differential form $\omega$.

Consider the exponential map $\exp:S\R^n\times\R\rightarrow\R^n$ given by $\exp(x,v,t):=x+tv$. Then it is easy to see that the pullback of the volume form on $\R^n$ along the exponential map satisfies
\begin{align*}
	\exp^*(dx_1\wedge\dots\wedge dx_n)=\left(\sum_{i=0}^{n-i}t^{n-i-1}\tilde{\kappa}_i\right)\wedge (\alpha+dt).
\end{align*}
for certain $\SO(n)$-invariant differential forms $\tilde{\kappa}_i\in \Omega^{n-1}(S\R^n)$, compare \cite{FuAlgebraicIntegralGeometry2014} Section 2.1. More precisely,
\begin{align*}
	\tilde{\kappa}_i=\frac{1}{i!(n-1-i)!}\sum_{\sigma\in S_n}\sign(\sigma)v_{\sigma(1)}dx_{\sigma(2)}\dots dx_{\sigma(i+1)}\wedge dv_{\sigma(i+2)}\dots dv_{\sigma(n)},
\end{align*}
where $(x_1,\dots,x_n)$ are coordinates on $\R^n$ with induced coordinates $(v_1,\dots,v_n)$ on $S^{n-1}\subset\R^n$ and $S_n$ denotes the symmetric group of degree $n$.
Then for $0\le i\le n-1$, the $i$th intrinsic volume $\mu_i$ is given by
\begin{align*}
	\mu_i(K)=\frac{1}{(n-i)\omega_{n-i}}\nc(K)[\tilde{\kappa}_i],
\end{align*}
where $\omega_k$ denotes the $k$-dimensional volume of the unit ball in $\R^k$.
Using this normalization, $\mu_i(B_1(0))=\binom{n}{i}\frac{\omega_n}{\omega_{n-i}}$.
It is easy to see that the curvature measure induced by $\tilde{\kappa}_i$ is non-negative. In fact, these curvature measures are precisely the Federer curvature measures \cite{FedererCurvatureMeasures1959}. As the intrinsic volumes are monotone, this is also a consequence of Lemma \ref{lemma:monotoneValuationsCurvatureMeasures} and the following simple relation, which follows from an elementary calculation (see also \cite{FuSomeremarksLegendrian1998} Lemma 3.1).
\begin{lemma}
	\label{lemma:intrinsicVolumesDifferentialForms}
	\begin{align*}
		d\tilde{\kappa}_i&=(n-i)\alpha\wedge \tilde{\kappa}_{i-1}\quad \text{for }1\le i\le n-1,\\
		d\tilde{\kappa}_0&=0.
	\end{align*}
	In particular, for $1\le i\le n-1$,
	\begin{equation}
		\label{eq:RuminDiffTildeKappa}
		i_TD\tilde{\kappa}_i=(n-i)\tilde{\kappa}_{i-1}.
	\end{equation}
\end{lemma}

\subsection{Convex functions}	
\label{section:convexFunctions}
We refer to the monographs by Rockafellar \cite{RockafellarConvexanalysis1997} and Rockafellar and Wets \cite{RockafellarWetsVariationalanalysis1998} for a general background on convex functions and only collect some basic facts. We will only be interested in finite-valued convex functions. These are in particular locally Lipschitz continuous, as the following special case of \cite{RockafellarWetsVariationalanalysis1998} 9.14 shows.
\begin{proposition}
	\label{proposition:bound_lipschitz_constant}
	Let $U\subset \R^n$ be a convex open subset and $f:U\rightarrow\R$ a convex function. If $X\subset U$ is a set with $X+\epsilon B_1(0)\subset U$ and $f$ is bounded on $X+ \epsilon B_1(0)$, then $f$ is Lipschitz continuous on $X$ with Lipschitz constant $\frac{2}{\epsilon}\sup_{x\in X+\epsilon B_1(0)}|f(x)|$.
\end{proposition}

The space $\Conv(\R^n,\R)$ of finite-valued convex functions on $\R^n$ is naturally equipped with the metrizable topology induced by epi-convergence. For our purposes, the precise definition of epi-convergence is not relevant, as this type of convergence has a very simple description on $\Conv(\R^n,\R)$: A sequence $(f_j)_j$ in $\Conv(\R^n,\R)$ epi-converges to $f\in\Conv(\R^n,\R)$ if and only if it converges uniformly to $f$ on compact subsets of $\R^n$. This in turn is the case if and only if $(f_j)_j$ converges pointwise to $f$, see \cite{RockafellarWetsVariationalanalysis1998} Theorem 7.17. This implies in particular that the map $f\mapsto  \sup_{x\in A}|f(x)|$ is continuous on $\Conv(\R^n,\R)$ for all compact subsets $A\subset\R^n$.\\

For $f\in\Conv(\R^n,\R)$, its convex conjugate or Legendre transform $\mathcal{L}f:\R^n\rightarrow(-\infty,\infty]$ is defined by
\begin{align*}
	\mathcal{L}f(y)=\sup_{x\in\R^n}\left(\langle y,x\rangle -f(x)\right).
\end{align*}
If $f$ is of class $C^\infty$ and such that there exists a constant $\lambda>0$ such that $x\mapsto f(x)-\lambda\frac{|x|^2}{2}$ is convex, then $df:\R^n\rightarrow\R^n$ is a diffeomorphism and
\begin{align*}
	\mathcal{L}f(df(x))=\langle df(x),x\rangle -f(x) \quad\text{for all }x\in\R^n.
\end{align*} 
Using the inverse function theorem, it is easy to see that this also implies that $\mathcal{L}f$ is of class $C^\infty$ with positive definite Hessian.\\

Let us consider the following two subspaces of $\Conv(\R^n,\R)$:
\begin{align*}
	\Conv_0(\R^n,\R):=&\{f\in\Conv(\R^n,\R): f(0)=0< f(x)\quad\text{for all } 0\ne x\in\R^n\},\\
	\Conv_0^+(\R^n,\R):=&\left\{f\in\Conv_0(\R^n,\R):f-\lambda\frac{|\cdot|^2}{2}\in\Conv(\R^n,\R) \text{ for some }\lambda>0\right\}.
\end{align*}
Note that $f\in\Conv_0^+(\R^n,\R)$ implies that the sublevel sets
$\{y\in\R^n: \mathcal{L}f(y)\le t\}$ are compact and non-empty for $t\ge 0$. Moreover, if $f\in \Conv_0^+(\R^n,\R)$ is of class $C^1$, then $df(x)=0$ if and only if $x=0$, as $f$ has a unique minimum in $x=0$. If $f$ is in addition smooth, this also implies $d\mathcal{L}f(y)=0$ if and only if $y=0$.\\

We will need the following estimates in the proof of Theorem \ref{theorem:existenceFunctionalIntrinsicVolumes}.
\begin{lemma}
	\label{lemma:bounds_subgradient}
	Let $f\in\Conv^+_0(\R^n,\R)$ be a smooth function. If $0\ne y=df(x_0)$ for $x_0\in\R^n$, then
	\begin{enumerate}
		\item $|\mathcal{L}f(y)|\le (1+2|x_0|)\sup_{|x|\le |x_0|+1}|f(x)|$,
		\item $|y|\le 2\sup_{|x|\le |x_0|+1}|f(x)|$.
	\end{enumerate}
	Moreover, $\mathcal{L}f(y)\le (1+2R)\sup_{|x|\le R+1}|f(x)|$ implies $|y|\le 2\sup_{|x|\le R+1}|f(x)|$ for $R>0$.
\end{lemma}
\begin{proof}
	As $y=df(x_0)$, we may bound $|y|$ by the Lipschitz constant of $f$ on a neighborhood of $x_0$, which is bounded by $2\sup_{|x|\le |x_0|+1}|f(x)|$ by Proposition \ref{proposition:bound_lipschitz_constant}. In addition,
	\begin{align*}
		\mathcal{L}f(y)+f(x_0)=\langle y,x_0\rangle 
	\end{align*}
	holds, so the first estimate follows from
	\begin{align*}
		|\mathcal{L}f(y)|\le |f(x_0)|+|x_0|\cdot|y|\le (1+2|x_0|)\sup_{|x|\le |x_0|+1}|f(x)|.
	\end{align*}
	For the second estimate, assume that $0\ne y\in \R^n$ is a point with $\mathcal{L}f(y)\le (1+2R)\sup_{|x|\le R+1}|f(x)|$. Then
	\begin{align*}
		\langle y,x_1\rangle -f(x_1)\le \mathcal{L}f(y)\le (1+2R)\sup_{|x|\le R+1}|f(x)|
	\end{align*}
	for all $x_1\in {\R^n}$. We may thus choose $x_1=(R+1) \frac{y}{|y|}$ to obtain
	\begin{align*}
		(R+1)|y|-\sup_{|x|\le R+1}f(x)\le (1+2R)\sup_{|x|\le R+1}|f(x)|.
	\end{align*}
	As $f$ is non-negative, this implies
	\begin{align*}
		(R+1)|y|\le (2+2R)\sup_{|x|\le R+1}|f(x)|,
	\end{align*}
	which shows $|y|\le 2\sup_{|x|\le R+1}|f(x)|$.
\end{proof}

\subsection{The classes $D^n_i$}
\label{section:Dnk}
Recall that we have defined the norm $\|\zeta\|$ for $\zeta\in D^n_i$, $1\le i\le n-1$, by
\begin{align*}
	\|\zeta\|=\sup_{t>0}\left|(n-i)\int_t^\infty \zeta(s)s^{n-i-1}ds\right|+\sup_{t>0}\left|t^{n-i}\zeta(t)+(n-i)\int_t^\infty \zeta(s)s^{n-i-1}ds\right|.
\end{align*}
For $R>0$ we also use $D^n_{i,R}$ for the subspaces of functions that are supported on $(0,R]$ and we will consider $C_c([0,\infty))$ as a subspace of $D^n_i$ using the obvious inclusion.
\begin{lemma}
	\label{lemma:density_cont_Dnk}
	For $\zeta\in D^{n}_{i,R}$ and $r>0$ define $\zeta^r\in C_c([0,\infty))\cap D^n_{i,R}$ by
	\begin{align*}
		\zeta^r(t):=\begin{cases}
			\zeta(t) & \text{for }t> r,\\
			\zeta(r) & \text{for }0\le t\le r.
		\end{cases}
	\end{align*}
	Then $\lim\limits_{r\rightarrow0}\|\zeta-\zeta^r\|=0$.
\end{lemma}
\begin{proof}
	For $t>0$ set
	\begin{align*}
		\rho_r(t):=&(n-i)\int_t^\infty s^{n-i-1}\zeta^r(s)ds, &&\eta_r(t):=t^{n-i}\zeta^r(t)+\rho_r(t),\\
		\rho(t):=&(n-i)\int_t^\infty s^{n-i-1}\zeta(s)ds, &&\eta(t):=t^{n-i}\zeta(t)+\rho(t).
	\end{align*}
	Then $\rho$ and $\eta$ extend to continuous functions on $[0,\infty)$ by assumption. Thus
	\begin{align*}
		\|\zeta-\zeta^r\|=&\|\rho-\rho_r\|_\infty+\|\eta-\eta_r\|_\infty
		\le 2\sup_{t\in[0,\infty)}|\rho(t)-\rho_r(t)|+\sup_{t\in[0,\infty)}|t^{n-i}\zeta(t)-t^{n-i}\zeta^r(t)|\\
		=&2\sup_{t\in [0,r]}|\rho(t)-\rho_r(t)|+\sup_{t\in [0,r]}|t^{n-i}\zeta(t)-t^{n-i}\zeta(r)|.
	\end{align*} 
	The second term obviously converges to $0$. For the first term, note that for $t\in[0,r]$,
	\begin{align*}
		\rho(t)-\rho_r(t)=&(n-i)\left[\int_t^r s^{n-i-1}\zeta(s)ds-\int_t^r s^{n-i-1}\zeta^r(s)ds\right]\\
		=&\rho(t)-\rho(r)-(n-i)\zeta(r)\int_t^rs^{n-i-1}ds\\
		=&\rho(t)-\rho(r)-\zeta(r)(r^{n-i}-t^{n-i}).
	\end{align*}
	As $\rho$ is continuous in $0$, $\sup_{t\in[0,r]}|\rho(t)-\rho(r)|$ converges to $0$ for $r\rightarrow0$. As $\zeta\in D^n_i$,  $\sup_{t\in[0,r]}|\zeta(r)|(r^{n-i}-t^{n-i})\le 2r^{n-i}|\zeta(r)|$ also converges to $0$. The claim follows. 
\end{proof}

\begin{lemma}
	$D^n_{i,R}$ is a Banach space with respect to $\|\cdot\|$. Moreover, $D^n_{i,R}$ is the completion of $C_c([0,\infty))\cap D^n_{i,R}$ with respect to $\|\cdot\|$.
\end{lemma}
\begin{proof}
	Let $(\zeta_j)_j$ be a Cauchy sequence in $D^n_{i,R}$. Then for $t_0>0$
	\begin{align*}
		&t_0^{n-i}|\zeta_j(t_0)-\zeta_l(t_0)|\le 2\|\zeta_j-\zeta_l\|.
	\end{align*}
	Consider the functions $\phi_j(t):=t^{n-i}\zeta_j(t)$, $t>0$. As $\zeta_j\in D^n_i$, $\phi_j$ extends to an element of $C_c([0,\infty))$ by $\phi_j(0):=0$. Moreover, $\|\phi_j-\phi_l\|_\infty\le 2\|\zeta_j-\zeta_l\|$, so this is a Cauchy sequence in $C_c([0,\infty))$. As all functions are supported on $[0,R]$, this sequence converges uniformly to some $\phi\in C_c([0,\infty))$ that is supported on $[0,R]$ and satisfies $\phi(0)=\lim\limits_{j\rightarrow\infty}\phi_j(0)=0$.\\
	
	In particular, $(\zeta_j)_j$ converges locally uniformly on $(0,\infty)$ to the function $\zeta\in C_b((0,\infty))$ given by $\zeta(t)=\frac{\phi(t)}{t^{n-i}}$ for $t>0$. This shows that
	\begin{align*}
		\int_t^\infty\zeta(s)s^{n-i-1}ds&=\lim\limits_{j\rightarrow\infty}\int_t^\infty \zeta_j(s)s^{n-i-1}ds \quad\text{for } t>0.
	\end{align*}
	From the definition of the norm $\|\cdot\|$, we thus obtain $\|\zeta\|\le \liminf_{l\rightarrow\infty}\|\zeta_l\|$, or more generally
	\begin{align*}
		\|\zeta-\zeta_j\|\le \liminf_{l\rightarrow\infty}\|\zeta_l-\zeta_j\|.
	\end{align*}
	This implies that $(\zeta_j)_j$ converges to $\zeta$ with respect to $\|\cdot\|$: Let $\epsilon>0$ be given. As $(\zeta_j)_j$ is a Cauchy sequence with respect to $\|\cdot\|$, there exists $j_0>0$ such that $\|\zeta_j-\zeta_l\|\le \epsilon $ for $j,l\ge j_0$. Thus $j\ge j_0$ implies
	\begin{align*}
		\|\zeta-\zeta_j\|\le&\liminf\limits_{l\rightarrow\infty}\|\zeta_l-\zeta_j\|\le \epsilon.
	\end{align*}
	Thus the sequence converges to $\zeta$.\\
	
	It remains to see that $\zeta\in D^n_i$. First note that
	\begin{align*}
		\lim\limits_{t\rightarrow0}t^{n-i}\zeta(t)&=\phi(0)=0.
	\end{align*}
	By definition of the norm $\|\cdot\|$, the continuous functions $\psi_j\in C_c([0,\infty))$,
	\begin{align*}
		\psi_j(t):=\int_t^\infty \zeta_j(s)s^{n-i-1}ds,
	\end{align*}
	converge uniformly on $[0,\infty)$ to $\psi(t)=\int_t^\infty \zeta(s)s^{n-i-1}ds$, which is thus also continuous. Thus the limit $\psi(0)=\lim\limits_{t\rightarrow0}\int_t^\infty \zeta(s)s^{n-i-1}ds$ exists and is finite, which shows $\zeta\in D^n_i$.\\
	
	The last claim follows from the fact that $C_c([0,\infty))\cap D^n_{i,R}$ is dense in $D^n_{i,R}$ by Lemma \ref{lemma:density_cont_Dnk}.
\end{proof}

\section{Dually epi-translation invariant valuations on convex functions}
\label{section:VConv}
\subsection{Topology and support}
\label{section:duallyTopologySupport}
The existence of the homogeneous decomposition $
\VConv(\R^n)=\bigoplus\limits_{i=0}^{n}\VConv_i(\R^n)$
implies that the polarization of $\mu\in\VConv_i(\R^n)$ given by
\begin{align*}
	\bar{\mu}(f_1,\dots,f_i)=\frac{1}{i!}\frac{\partial^i}{\partial \lambda_1\dots\partial \lambda_i}\Big|_0\mu\left(\sum_{i=1}^{i}\lambda_if_i\right)\quad \text{for }f_1,\dots, f_i\in\Conv(\R^n,\R)
\end{align*}
is well defined, as the right hand side of this equation is a polynomial in $\lambda_1,\dots,\lambda_i\ge 0$.

One of the most important constructions derived from the polarization is the construction of certain distributions associated to homogeneous valuations (for translation invariant valuations on convex bodies, this construction is due to Goodey and Weil, see \cite{GoodeyWeilDistributionsvaluations1984}).
\begin{theorem}[\cite{Knoerrsupportduallyepi2021} Theorem 2]
	\label{theorem:Def_GW}
	For every $\mu\in\VConv_i(\R^n)$ there exists a unique distribution $\GW(\mu)\in\mathcal{D}'((\R^n)^i)$ with compact support which satisfies the following property: If $f_1,...,f_i\in \Conv(\R^n,\R)\cap C^\infty(\R^n)$, then
	\begin{align*}
		\GW(\mu)[f_1\otimes...\otimes f_i]=\bar{\mu}(f_1,...,f_i).
	\end{align*}
	Moreover, the support of $\GW(\mu)$ is contained in the diagonal in $(\R^n)^i$.
\end{theorem}
Here, $f_1\otimes\dots\otimes f_i\in C^\infty((\R^n)^i)$ is given by $[f_1\otimes\dots\otimes f_i](x_1,\dots,x_i)=f_1(x_1)\dots f_i(x_i)$ for $x_1,\dots,x_i\in\R^n$.
The distribution $\GW(\mu)$ is called the Goodey-Weil distribution of $\mu\in\VConv_i(\R^n)$. We thus obtain a well defined and injective map $\GW:\VConv_i(\R^n)\rightarrow\mathcal{D}'((\R ^n)^i)$, called Goodey-Weil embedding.
Moreover, $\GW(\mu)$ satisfies the inequality (compare Section 5.1 in \cite{Knoerrsupportduallyepi2021})
\begin{align}
	\label{inequality:GW_continuity}
	|\GW(\mu)[\phi_1\otimes\dots\otimes\phi_i]|\le c_i \|\mu\|_C\prod_{k=1}^i\|\phi_k\|_{C^2(\R^n)}
\end{align}
for $\phi_1,\dots,\phi_i\in C^\infty_c(\R^n)$, where $c_i>0$ is a constant depending on $i$ only, $C\subset \Conv(\R^n,\R)$ is a certain compact subset (see \cite{Knoerrsupportduallyepi2021} Proposition 2.4 for a description of these subsets) and 
\begin{align}
	\label{eq:defSeminorm}
	\|\mu\|_C:=\sup_{f\in C}|\mu(f)|
\end{align}
defines a continuous norm on $\VConv_i(\R^n)$. This immediately implies 
\begin{corollary}
	\label{corollary:weak_continuity_GW}
	For $\phi_1,\dots,\phi_i\in C^\infty_c(\R^n)$ the map
	\begin{align*}
		\VConv_i(\R^n)&\rightarrow\R\\
		\mu&\mapsto \GW(\mu)[\phi_1\otimes\dots\otimes\phi_i]
	\end{align*}
	is continuous.
\end{corollary}

In \cite{Knoerrsupportduallyepi2021} the support of $\mu\in\VConv(\R^n)$ was defined as $\supp\mu:=\bigcup_{i=1}^n\Delta_i^{-1}(\supp\GW(\mu_i))$, where $\mu=\sum_{i=0}^{n}\mu_i$ is the decomposition of $\mu$ into its homogeneous components and $\Delta_i:\R^n\rightarrow (\R ^n)^i$ is the diagonal embedding. In particular, $\supp\mu=\emptyset$ for $\mu\in\VConv_0(\R^n)$ by definition. Note that $\supp\mu$ is a compact subset, as the Goodey-Weil distributions have compact support. The fact that the support of these distributions is contained in the diagonal leads to the following alternative characterization of the support:
\begin{proposition}[\cite{Knoerrsupportduallyepi2021} Proposition 6.3]
	\label{proposition:supportValuation}
	The support of $\mu\in\VConv(\R^n)$ is minimal (with respect to inclusion) among the closed sets $A\subset \R^n$ with the following property: If $f,h\in\Conv(\R^n,\R)$ satisfy $f=h$ on an open neighborhood of $A$, then $\mu(f)=\mu (h)$.
\end{proposition}
For a closed subset $A\subset \R^n$, we will denote by $\VConv_{A}(\R^n)\subset\VConv(\R^n)$ the subspace consisting of valuations with support contained in $A$.\\
Recall that $\VConv(\R^n)$ and its subspaces are equipped with the topology of uniform convergence on compact subsets in $\Conv(\R^n,\R)$. This topology is induced by the semi-norms $\|\cdot\|_C$ defined in \eqref{eq:defSeminorm} for all $C\subset\Conv(\R^n,\R)$ compact. For spaces of valuations with compact support, the subspace topology is actually much simpler:
\begin{theorem}[\cite{Knoerrsupportduallyepi2021} Proposition 6.8]
	\label{thm:BanachStructureCompact}
	If $A\subset \R^n$ is compact, then $\VConv_A(\R^n)$ is
	a Banach space.
	If $A$ is in addition convex, then the subspace topology is induced by the continuous norm
	\begin{align*}
		\|\mu\|_{A,1}:=\sup\left\{|\mu(f)|  \ : \ f\in\Conv(\R^n,\R), \sup\limits_{\mathrm{dist}(x,A)\le 1}|f(x)|\le 1 \right\}.
	\end{align*}
\end{theorem}
Note that this norm is only well-defined for valuations with support in $A$ (or more generally support contained in $\{x\in\R^n:\mathrm{dist}(x,A)<1\}$).
When considering convergent sequences in $\VConv(\R^n)$, it is usually  sufficient to consider spaces of valuations with certain bounds on the support, as shown by the following proposition.
\begin{proposition}[\cite{Knoerrsupportduallyepi2021} Proposition 1.2]
	\label{proposition:restriction_support_convergent_sequences}
	If a sequence $(\mu_j)_j$ in $\VConv(\R^n)$ converges to $\mu$ in $\VConv(\R^n)$, then there exists a compact set $A\subset \R^n$ such that the supports of $\mu$ and $\mu_j$ are contained in $A$ for all $j\in\mathbb{N}$. In particular, $(\mu_j)_j$ converges to $\mu$ in $\VConv_A(\R^n)$.
\end{proposition}
The following observation will play a key role in the proof of Theorem \ref{maintheorem:continuousHadwiger}.
\begin{proposition}[\cite{Knoerrsupportduallyepi2021} Proposition 6.4]
	\label{proposition:shape_support}
	If the support of $\mu\in  \VConv(\R^n)$ is contained in a one-point set, then it is empty and $\mu$ is constant.
\end{proposition}
For non-trivial $i$-homogeneous valuations, one can show finer restrictions on the minimal dimension of the support, see \cite{KnoerrUnitarilyinvariantvaluations2021} Section 2.2, however, we will only need this coarser result.
\subsection{Smooth valuations}
\label{section:smoothValuations}
For a smooth convex function $f\in\Conv(\R^n,\R)$, the graph  of its differential defines a smooth oriented submanifold of the cotangent bundle $T^*\R^n\cong\R^n\times\R^n$. We can consider this submanifold as an $n$-current by integrating suitable differential forms over it. For arbitrary convex functions, this graph exists in a generalized sense as shown by Fu.
\begin{theorem}[\cite{FuMongeAmperefunctions.1989} Theorem 2.0 and Proposition 3.1]
	For every $f\in\Conv(\R^n,\R)$ there exists a unique integral current $D(f)\in I_n(T^*\R^n)$ with the following properties:
	\begin{enumerate}
		\item $D(f)$ is closed, i.e. $\partial D(f)=0$,
		\item $D(f)$ is Lagrangian, i.e. $D(f)$ vanishes on multiples of $\omega_s$,
		\item $D(f)$ is locally vertically bounded, i.e. $\supp D(f)\cap \pi^{-1}(A)$ is compact for all $A\subset \R^n$ compact,
		\item $D(f)(\phi(x,y)\pi^*\vol)=\int_{\R^n}\phi(x,df(x))d\vol(x)$ for all $\phi\in C^\infty_c(T^*\R^n)$.
	\end{enumerate}
\end{theorem}
Here $\pi:T^*\R^n\cong\R^n\times\R^n\rightarrow\R^n$ denotes the projection onto the first factor, $\vol$ is the natural volume form on $\R^n$, $\omega_s$ denotes the natural symplectic form on $\R^n\times\R^n$, and the boundary $\partial S$ of a current $S$ is the current defined by $\partial S(\omega)=S(d\omega)$ for compactly supported forms $\omega$. Moreover, $I_n(T^*\R^n)$ denotes the space of integral $n$-currents on $T^*\R^n$. We refer to \cite{FedererGeometricmeasuretheory1969} for the precise definition of this notion and for a general background on currents. The current $D(f)$ is called the differential cycle of $f$ and is given by integration over the graph of $df$ if $f$ is smooth.\\
More generally, a function is called \emph{Monge-Amp\`ere} if a current with these properties exists. This class of functions is much larger than the class of convex functions and we refer to \cite{FuMongeAmperefunctions.1989,JerrardSomerigidityresults2010} for a detailed study of this notion.\\ 

As shown in \cite{KnoerrSmoothvaluationsconvex2024} Section 4, every differential form $\omega\in \Omega^{n-i}_c(\R^n)\otimes\Lambda^{i}\R^n\subset\Omega^n(\R^n\times\R^n)$ defines a continuous, dually epi-translation invariant valuation $\mu\in\VConv_i(\R^n)$ by $	\mu(f):=D(f)[\omega]$. Here $\Omega_c^{n-i}(\R^n)$ denotes the space of compactly supported $(n-i)$-forms on $\R^n$ and $\Lambda^i\R^n$ the space of differential $i$-forms on $\R^n$ with constant coefficients. Valuations of this type are called smooth. By the main result of \cite{KnoerrSmoothvaluationsconvex2024}, smooth valuations form a dense subspace of $\VConv_i(\R^n)$. We will need the following version of this result for invariant valuations:
\begin{theorem}[\cite{KnoerrSmoothvaluationsconvex2024} Proposition 6.6]
	\label{theorem:approximation_smooth_valuations}
	Let $G\subset \SO(n)$ be a compact subgroup. Then the space of smooth $G$-invariant valuations is sequentially dense in the space of all $G$-invariant valuations in $\VConv_i(\R^n)$.
\end{theorem}
In fact, every such valuation may be represented by a smooth $G$-invariant differential form, see \cite{KnoerrSmoothvaluationsconvex2024} Corollary 6.7. Furthermore, the Lefschetz decomposition for the symplectic manifold $T^*\R^n$ implies that one may also assume that such a differential $n$-form $\tau$ is \emph{primitive}, that is, $\tau\wedge \omega_s=0$.\\

The differential form inducing a given smooth valuation is highly non-unique. The kernel of this procedure can be described using a certain second order differential operator $\D$ naturally defined on any symplectic manifold, which we will call the \emph{symplectic Rumin differential}. It is also one of the operators considered in \cite{TsengYauCohomologyHodgetheory2012}.
\begin{theorem}[\cite{KnoerrSmoothvaluationsconvex2024} Theorem 2]
	\label{theorem:kernel}
	$\tau\in\Omega^{n-i}_c(\R^n)\otimes\Lambda^{i}\R^n$ satisfies $D(f)[\tau]=0$ for all $f\in\Conv(\R^n,\R)$ if and only if
	\begin{enumerate}
		\item $\D\tau=0$,
		\item $\int_{\R^n}\tau=0$, where we consider the zero section $\R^n\hookrightarrow T^*\R^n$ as a submanifold.
	\end{enumerate}
\end{theorem}
Note that the first condition is always satisfied for $i=0$, while the second is always satisfied for $1\le i\le n$. We will not need the precise definition of the operator $\D$, only that it vanishes on closed forms as well as multiples of the symplectic form $\omega_s$, see \cite{KnoerrSmoothvaluationsconvex2024} Proposition 4.13. 
\section{Properties of the Differential cycle}
\label{section:PropDifferentialCycle}
\subsection{Weak continuity of the differential cycle on $\Conv(\R^n,\R)$}
In this section, we collect some results on the differential cycle and show that it satisfies a weak continuity property. These results are strictly speaking not necessary for the construction of the singular valuations in Section \ref{section:singularHessianValuations}, as our construction only uses that the Hessian measures (which can be obtained by integration over the differential cycle, see Proposition \ref{proposition:interpretation-hessian-measures} below) have similar properties, compare for example \cite{ColesantiEtAlHessianvaluations2020}. However, a generalization of this construction to different invariant measures will need these continuity properties of the differential cycle for arbitrary continuous differential forms. We have thus included the argument.\\

First, we will collect some properties of the differential cycle, established by Fu in \cite{FuMongeAmperefunctions.1989}.
\begin{proposition}[\cite{FuMongeAmperefunctions.1989} Proposition 2.4]
	\label{proposition_Fu_sum_MA_andC11}
	Let $f:\R^n\rightarrow \R$ be a Monge-Amp\`ere function and $\phi\in C^{1,1}(\R^n)$. Then $f+\phi$ is Monge-Amp\`ere and 
	\begin{align*}
		D(f+\phi)=G_{\phi*}D(f),
	\end{align*}
	where $G_\phi:T^*\R^n\rightarrow T^*\R^n$ is given by $(x,y)\mapsto (x,y+d\phi(x))$.
\end{proposition}
\begin{proposition}[\cite{FuMongeAmperefunctions.1989} Theorem 2.2.]
	\label{proposition:differential_cycle_estimate_support}
	If $f\in\Conv(\R^n,\R)$, then 
	\begin{align*}
		\supp D(f)\subset \mathrm{graph}\ \partial f:=\left\{(x,y)\in T^*\R^n: y\in\partial f(x)\right\},
	\end{align*}
	where $\partial f(x)$ denotes the subgradient of $f$ in $x\in {\R^n}$. In particular, given an open set $U\subset {\R^n}$,
	\begin{align*}
		\supp D(f)\cap\pi^{-1}(U)\subset U\times B_{\mathrm{lip}(f|_U)}(0),
	\end{align*}
	where $\mathrm{lip}(f|_U)$ denotes the Lipschitz constant of $f|_U$.
\end{proposition}
For $R>0$ let $U_R(0)$ denote the open ball in $\R^n$ centered at the origin.
\begin{lemma}[\cite{KnoerrSmoothvaluationsconvex2024} Lemma 5.8]
	\label{lemma:differential_cycle_mass_estimate}
	For $f\in\Conv(\R^n,\R)$,
	\begin{align*}
		M_{\pi^{-1}(U_R(0))}(D(f))\le2^n\omega_n\left(\sup\limits_{|x|\le R+1}|f(x)|\right)^n.
	\end{align*}
\end{lemma}
Here, 
\begin{align*}
	M_U(S):=\sup\{|S(\omega)|: \omega\in\Omega^n(T^*\R^n)\text{ with } \supp\omega\subset U, \|\omega\|_\infty\le 1\}
\end{align*}
denotes the mass of an $n$-current $S$ on an open subset $U\subset T^*\R^n$, where $\|\omega\|_\infty:=\sup_{(x,y)\in T^*\R^n} \left\|\omega|_{(x,y)}\right\|$.\\
	
	We will need the following continuity property of the differential cycle.
	\begin{theorem}[\cite{KnoerrSmoothvaluationsconvex2024} Theorem 5.9]
		\label{theorem:differential_cycle_continuity_local_flat_topology}
		$D:\Conv(\R^n,\R)\rightarrow I^n(T^*\R^n)$ is continuous with respect to the local flat topology on $I^n(T^*\R^n)$. 
	\end{theorem}
	We refer to \cite{FedererGeometricmeasuretheory1969,KnoerrSmoothvaluationsconvex2024} for the definition of the semi-norms defining the local flat topology and only note that this implies that $f\mapsto D(f)[\phi(x,y) \omega]$ defines a continuous function for every $\omega\in \Omega^n(T^*\R^n)$ and $\phi\in C_c^1(T^*\R^n)$. We are now able to establish the desired continuity properties of the differential cycle.
	\begin{proposition}
		\label{proposition:vague_continuity_differential_cycle}
		Let $\omega\in\Omega^n(T^*\R^n)$. If $(f_j)_j$ is a sequence in $\Conv(\R^n,\R)$ converging to $f\in\Conv(\R^n,\R)$, then 
		\begin{align*}
			\lim\limits_{j\rightarrow\infty}D(f_j)[\phi \omega]=D(f)[\phi\omega]\quad\text{for all } \phi\in C_c(T^*\R^n).
		\end{align*}
	\end{proposition}
	\begin{proof}
		By multiplying $\omega$ with a smooth cut-off function, it is enough to show the claim for $\omega\in \Omega^n_c(T^*\R^n)$, so let us assume that $\supp\omega\subset U_R(0)\times U_R(0)\subset T^*\R^n$ for some $R>0$. Let $(f_j)_j$ be a sequence in $\Conv(\R^n,\R)$ that converges to $f\in\Conv(\R^n,\R)$. Then $(f_j)_j$ converges uniformly to $f$ on $B_{R+1}(0)$. In particular, there exists $C>0$ such that $\sup_{|x|\le R+1}|f_j(x)|\le C$ for all $j\in\mathbb{N}$. \\			
		Define a sequence of linear functionals $S_j:C(B_R(0)\times B_R(0))\rightarrow\R$ by $S_j(\phi):=D(f_j)[\phi\omega]$. Here we use that $\phi\omega$ is supported on $U_R(0)\times U_R(0)$ and we trivially extend this form to $T^*\R^n$. We also set $S(\phi):=D(f)[\phi\omega]$. Lemma \ref{lemma:differential_cycle_mass_estimate} thus implies 
		\begin{align*}
			|S_j(\phi)|=&|D(f_j)[\phi\omega]|\le M_{\pi^{-1}(U_R(0))}(D(f_j))\|\phi\|_\infty \|\omega\|_\infty\le 2^n\omega_n\left(\sup\limits_{|x|\le R+1}|f_j(x)|\right)^n\|\phi\|_\infty\|\omega\|_\infty\\
			\le&2^n\omega_nC^n\|\phi\|_\infty\|\omega\|_\infty.
		\end{align*}		
		In particular, the operator norm of $S_j$ is bounded by a constant independent of $j\in\mathbb{N}$.\\
		Next, consider $\phi\in C^1(B_R(0)\times B_R(0))$. Theorem \ref{theorem:differential_cycle_continuity_local_flat_topology} implies that in this case
		\begin{align*}
			\lim\limits_{j\rightarrow\infty}D(f_j)[\phi\omega]=D(f)[\phi\omega].
		\end{align*}
		Because $C^1(B_R(0)\times B_R(0))\subset C(B_R(0)\times B_R(0))$ is dense, the sequence $(S_j)_j$ converges pointwise on a dense subspace of $C(B_R(0)\times B_R(0))$. As the operator norms of these functionals are uniformly bounded, the sequence $(S_j)_j$ thus converges pointwise on $C(B_R(0)\times B_R(0))$ to some $\tilde{S}\in C(B_R(0)\times B_R(0))'$. But then $\tilde{S}=S$, as both sides are continuous and coincide on the dense subspace $C^1(B_R(0)\times B_R(0))$. Thus \begin{align*}
			\lim\limits_{j\rightarrow\infty}D(f_j)[\phi\omega]=\lim\limits_{j\rightarrow\infty}S_j(\phi)=S(\phi)=D(f)[\phi\omega]\quad\text{for all } \phi\in C(B_1(0)\times B_1(0)),
		\end{align*}
		which implies the claim.
	\end{proof}
	The restriction on the support of the functions can be relaxed to boundedness in the "horizontal" direction in the following sense: 
	\begin{theorem}
		\label{theorem:joint_continuity_differential_cycle}
		Let $(\phi_j)_j\subset C(T^*\R^n)$ be a sequence that converges locally uniformly to $\phi_0\in C(T^*\R^n)$ and let $(f_j)_j$ be a sequence in $\Conv(\R^n,\R)$ converging to $f_0\in\Conv(\R^n,\R)$. If there exists a compact subset $A\subset {\R^n}$ such that $\supp\phi_j\subset\pi^{-1}(A)$ for all $j\in\mathbb{N}$, then
		\begin{align*}
			\lim\limits_{j\rightarrow\infty}D(f_j)[\phi_j\omega]=D(f_0)[\phi_0\omega]
		\end{align*}
		for any $\omega\in \Omega^n(T^*\R^n)$.
	\end{theorem}
	\begin{proof}
		Choose $R>0$ such that  $A\subset U_R(0)$. Then the sequence $(f_j)_j$ converges uniformly to $f$ on $B_{R+1}(0)$, so Proposition \ref{proposition:bound_lipschitz_constant} shows that the Lipschitz constants of $f_j$ on $U_R(0)$ are bounded by some $L>0$ independent of $j\in\mathbb{N}$. Using Proposition \ref{proposition:differential_cycle_estimate_support}, we thus obtain $L>0$ such that  $\supp D(f_j)\cap \pi^{-1}(U_R(0))\subset U_R(0)\times B_L(0)$ for all $j\ge 0$. 
		Choose $\psi\in C_c(T^*\R^n,[0,1])$ with $\psi=1 $ on a neighborhood of $U_R(0)\times B_L(0)$. If $\phi\in C(\supp\psi)$, then
		\begin{align*}
			\tilde{\phi}(x,y):=\begin{cases}
				\psi(x,y)\phi(x,y)& \text{for }(x,y)\in \supp\psi,\\
				0 & \text{for }(x,y)\notin\supp\psi,
			\end{cases}
		\end{align*}
		defines a continuous function on $T^*\R^n$ with compact support as $\supp\psi$ is compact. Moreover, $\|\tilde{\phi}\|_\infty\le \|\phi\|_\infty$. We consider the map
		\begin{align*}
			\Conv(\R^n,\R)\times C(\supp\psi)&\rightarrow\R\\
			(f,\phi)&\mapsto D(f)[\tilde{\phi}\omega].
		\end{align*}
		By Proposition \ref{proposition:vague_continuity_differential_cycle} and the fact that $D(f)$ is an integral current, this is a separately continuous map. As $C(\supp\psi)$ is a Banach space and $\Conv(\R^n,\R)$ is metrizable, the principle of uniform boundedness implies that this map is jointly continuous. In particular,
		\begin{align*}
			\lim\limits_{j\rightarrow\infty}D(f_j)[\phi_j\omega]=	\lim\limits_{j\rightarrow\infty}D(f_j)[(\phi_j\psi)\omega]=D(f_0)((\phi_0\psi) \omega)=D(f_0)[\phi_0\omega].
		\end{align*}
	\end{proof}
	\begin{corollary}
		\label{corollary:differential_cycle_vague_continuity}
		Let $\omega\in \Omega^n(T^*\R^n)$. If $(f_j)_j$ is a sequence in $\Conv(\R^n,\R)$ converging to $f\in\Conv(\R^n,\R)$, then 
		\begin{align*}
			\lim\limits_{j\rightarrow\infty}D(f_j)[\phi \omega]=D(f)[\phi\omega]
		\end{align*}
		for all $\phi\in C(T^*\R^n)$ with $\pi(\supp\phi)\subset {\R^n}$ bounded.
	\end{corollary}
	\begin{proof}
		Apply Theorem \ref{theorem:joint_continuity_differential_cycle} to the constant sequence $\phi_j:=\phi$.
	\end{proof}	
	
	\subsection{The normal cycle of sublevel sets}	 
	\label{section:normalCyclesSublevel}
	For a convex body $K\in\mathcal{K}(\R^n)$ its support function $h_K\in\Conv(\R^n,\R)$ is defined by
	\begin{align*}
		h_K(y)=\sup_{x\in K}\langle y,x\rangle.
	\end{align*}
	If $K$ is a smooth convex body with strictly positive Gauss curvature, then $h_K$ is a smooth convex function on $\R^n\setminus\{0\}$ and if $\nu_K:\partial K\rightarrow S^{n-1}$ denotes the Gauss map, that is, the map that associates to $x\in \partial K$ its unique outer unit normal, then $h_K(y)=\langle y,\nu_K^{-1}(y)\rangle$.\\
	
	For $f\in \Conv_0^+(\R^n,\R)$ set $K_{\mathcal{L}f}^t:=\{y\in \R^n: \mathcal{L}f(y)\le t\}$ for $t>0$ and $h(\mathcal{L}f,x,t):=h_{K_{\mathcal{L}f}^t}(x)$ for $x\in {\R^n}$. The following is a version of \cite{ColesantiSalaniBrunnMinkowskiinequality2003} Theorem 4.
	\begin{proposition}
		\label{prop:differentialJointSupportFunctionSublevel}
		If $f\in\Conv_0^+(\R^n,\R)\cap C^\infty(\R^n)$, then $h(\mathcal{L}f,\cdot,\cdot)\in C^\infty(\R^n\setminus\{0\}\times(0,\infty))$ and 
		\begin{align*}
			\frac{1}{|d\mathcal{L}f(y)|}=\frac{\partial}{\partial t}h\left(\mathcal{L}f,\frac{d\mathcal{L}f(y)}{|d\mathcal{L}f(y)|},\mathcal{L}f(y)\right)\quad\text{for all } y\ne 0.
		\end{align*}
	\end{proposition}
	We recall the following notion, which will be used in the proof of the result below: Let $\Phi:M\rightarrow N$ be a smooth map between manifolds that is proper on the support of a current $S$ defined on $M$, that is, such that the intersection of the inverse image of any compact subset with $\supp S$ is compact. We then define a current $\Phi_\#S$ on $N$ by  $\Phi_\#S(\omega):=S(\Phi^*\omega)$ for $\omega\in \Omega_c^*(N)$.
	\begin{proposition}
		\label{prop:DerivativeNCSublevelsets}
		For $\omega\in \Omega^{n-1}(S\R^n)$, $f\in \Conv_0^+(\R^n,\R)\cap C^\infty(\R^n)$ and $t>0$
		\begin{align*}
			\frac{d}{ds}\Big|_0\nc(K_{\mathcal{L}f}^{t+s})[\omega]=\nc(K_{\mathcal{L}f}^t)\left[\frac{\partial}{\partial t}h(\mathcal{L}f,v,t)i_TD\omega\right]=\nc(K_{\mathcal{L}f}^t)\left[\frac{1}{|d\mathcal{L}f(y)|}i_TD\omega\right].
		\end{align*}
	\end{proposition}
	\begin{proof}	
		Note that $\mathcal{L}f\in C^\infty(\R^n)$ has positive definite Hessian, so the sublevel sets $K_{\mathcal{L}f}^t$ have $C^\infty$-boundary with strictly positive Gauss curvature for all $t>0$. In particular, for $t>0$
		\begin{align*}
			\frac{\partial}{\partial t}h(\mathcal{L}f,v,t)=\frac{\partial}{\partial t}h\left(\mathcal{L}f,\frac{d\mathcal{L}f(y)}{|d\mathcal{L}f(y)|},t\right)\quad \text{for }(y,v)\in\supp\nc(K_{\mathcal{L}f}^t),
		\end{align*}
		as $v=\frac{d\mathcal{L}f(y)}{|d\mathcal{L}f(y)|}$ is the unique outer unit normal to $K_{\mathcal{L}f}^t$ in $y\in\partial K_{\mathcal{L}f}^t$. Thus the second equation follows from Proposition \ref{prop:differentialJointSupportFunctionSublevel} once we have established the first equation.\\ 
		
		Let $\nu_{K_{\mathcal{L}f}^t}:\partial K_{\mathcal{L}f}^t\rightarrow S^{n-1}$ denote the Gauss map and $\chi(\cdot,t):\R^n\setminus\{0\}\rightarrow \R^n\setminus\{0\}$ the $0$-homogeneous extension of $\nu_{K_{\mathcal{L}f}^t}^{-1}$. Note that $\chi(\cdot,t)=\nabla_y h(\mathcal{L}f,\cdot,t)$. In particular, $\chi\in C^\infty(\R^n\setminus\{0\}\times(0,\infty))$ by Proposition \ref{prop:differentialJointSupportFunctionSublevel}.
		We define a smooth local flow $\Phi_s$ on $\R^{n}\setminus \{0\}$ for $s\in\R$ by 
		\begin{align*}
			\Phi_s(y):=\chi(d\mathcal{L}f(y),\mathcal{L}f(y)+s)\quad\text{for }y\in\R^{n}\setminus\{0\}\text{ with }\mathcal{L}f(y)\in (-s,\infty).
		\end{align*}
		By construction, $\Phi_s(K_{\mathcal{L}f}^t)=K_{\mathcal{L}f}^{t+s}$ for $s\in (-t,\infty)$. Note that $\Phi_s$ lifts to a local flow $\tilde{\Phi}_s$ on $(\R^n\setminus\{0\})\times S^{n-1}$ given for $(y,v)\in (\R^n\setminus\{0\})\times S^{n-1}$ by
		\begin{align*}
			\tilde{\Phi}_s(y,v):= \left(\chi(d\mathcal{L}f(y),\mathcal{L}f(y)+s),\frac{(d\Phi_s|_y^T)^{-1}(v)}{|(d\Phi_s|_y^T)^{-1}(v)|}\right)\quad \text{if }\mathcal{L}f(y)\in (-s,\infty),
		\end{align*}
		where $d\Phi_s|_y^T$ is the transpose of $d\Phi_s|_y$.	Then it is easy to see that $\tilde{\Phi}_s$ defines a diffeomorphism on a neighborhood of the support of $\nc(K_{\mathcal{L}f}^t)$ for all $s\in (-t,\infty)$. Moreover, the normal cycle of $\Phi_s(K_{\mathcal{L}f}^t)=K_{\mathcal{L}f}^{t+s}$ is the same as the image of the normal cycle of $K_{\mathcal{L}f}^t$ under $\tilde{\Phi}_s$, compare \cite{FuCurvaturemeasuressubanalytic1994} Remark 1.4.1. c).\\
		
		Let $\tilde{X}:=\frac{d}{ds}|_0\tilde{\Phi}_s$ and $X:=\frac{d}{ds}|_0\Phi_s$. For $(y,v)\in\supp\nc(K_{\mathcal{L}f}^t)$ we have $v=\frac{d\mathcal{L}f(y)}{|d\mathcal{L}f(y)|}$ and thus
		\begin{align*}
			\alpha(\tilde{X})|_{(y,v)}=&\langle v, X\rangle =\left\langle v,\frac{d}{ds}\Big|_0\chi (d\mathcal{L}f(y),t+s)\right\rangle =\frac{d}{ds}\Big|_0\langle v,\chi(d\mathcal{L}f(y),t+s)\rangle=\frac{d}{ds}\Big|_0h(\mathcal{L}f,v,t+s)\\
			=&\frac{\partial}{\partial t}h(\mathcal{L}f,v,t),
		\end{align*}
		where $\langle v,\chi(d\mathcal{L}f(y),t+s)\rangle=h(\mathcal{L}f,v,t+s)$ because $\chi(d\mathcal{L}f(y),t+s)=\chi(v,t+s)$ is the unique boundary point in $K_{\mathcal{L}f}^{t+s}$ such that $v$ is an outer unit normal to $K_{\mathcal{L}f}^{t+s}$.
		Without loss of generality we may assume that $d\omega=D\omega$. Thus
		\begin{align*}
			\frac{d}{ds}\Big|_0\nc(K_{\mathcal{L}f}^{t+s})[\omega]=&\frac{d}{ds}\Big|_0(\tilde{\Phi}_s)_{\#}\nc(K_{\mathcal{L}f}^t)[\omega]=\frac{d}{ds}\Big|_0\nc(K_{\mathcal{L}f}^t)[(\tilde{\Phi}_s)^*\omega]=\nc(K_{\mathcal{L}f}^t)[\mathcal{L}_{\tilde{X}}\omega]\\
			=&\nc(K_{\mathcal{L}f}^t)[i_{\tilde{X}}D\omega]=\nc(K_{\mathcal{L}f}^t)[i_{\tilde{X}}(\alpha\wedge i_TD\omega)]=\nc(K_{\mathcal{L}f}^t)[\alpha(\tilde{X})i_TD\omega]\\
			=&\nc(K_{\mathcal{L}f}^t)\left[\frac{\partial}{\partial t}h(\mathcal{L}f,v,t)i_TD\omega\right],
		\end{align*}
		where we have used that $D\omega=\alpha\wedge i_TD\omega$ as $D\omega$ is vertical and that the normal cycle of a convex body vanishes on closed forms as well as multiples of $\alpha$.
	\end{proof}
	Let $W^{1,2}((0,\infty))$ denote the Sobolev space of all functions in $L^2((0,\infty))$ with distributional derivatives in $L^2((0,\infty))$, $W_0^{1,2}((0,\infty))\subset W^{1,2}((0,\infty))$ the closure of $C^\infty_c((0,\infty))$ with respect to the norm $\|\phi\|_{W^{1,2}}:=\|\phi\|_{L^2}+\|\phi'\|_{L^2}$. 
	\begin{corollary}
		\label{corollary:inequality_partial_integration}
		Let $f\in\Conv_0^+(\R^n,\R)\cap C^\infty(\R^n)$. For $\phi\in W^{1,2}_0((0,\infty))$ and $c>0$
		\begin{align*}
			\int_0^c\phi(t)\nc(K_{\mathcal{L}f}^t)\left[\frac{1}{|d\mathcal{L}f(y)|}i_TD\omega\right]dt=\phi(c)\nc(K_{\mathcal{L}f}^c)[\omega]-\int_0^c \phi'(t)\nc(K_{\mathcal{L}f}^t)[\omega]dt.
		\end{align*}
		In particular, 
		\begin{align*}
			\left|\int_0^c	|\phi(t)|\nc(K_{\mathcal{L}f}^t)\left[\frac{1}{|d\mathcal{L}f(y)|}i_TD\omega\right]dt\right|\le2\int_0^c|\phi'(t)|dt\sup_{t\in [0,c]}|\nc(K_{\mathcal{L}f}^t)[w]|.
		\end{align*}
	\end{corollary}
	\begin{proof}
		By the Sobolev embedding, we have a continuous inclusion $W_0^{1,2}((0,\infty))\rightarrow C([0,\infty))$. Moreover, $\phi(t)=\int_0^t\phi'(s)ds$.
		The first equation now follows from Proposition \ref{prop:DerivativeNCSublevelsets} by partial integration for $\phi\in C_c^\infty((0,\infty))$. We may then pass to the closure of $C_c^\infty((0,\infty))$ in $W^{1,2}((0,\infty))$, as both sides depend continuously on $\phi\in W^{1,2}((0,\infty))$. This establishes the first equation.\\
		
		It is well known that $|\cdot|:W_0^{1,2}((0,\infty))\rightarrow W_0^{1,2}((0,\infty))$ is well defined and continuous with $|\phi|'=\phi'\cdot 1_{\{\phi>0\}}-\phi'\cdot 1_{\{\phi<0\}}$ almost everywhere (compare \cite{StampacchiaEquationselliptiquesdu1966}). Thus
		\begin{align*}
			&\left|\int_0^c|\phi(t)|\nc(K_{\mathcal{L}f}^t)[\frac{1}{|d\mathcal{L}f(y)|}i_TD\omega]dt\right|=\left||\phi(c)|\nc(K_{\mathcal{L}f}^c)-\int_0^c |\phi|'(t)\nc(K_{\mathcal{L}f}^t)[\omega]dt\right|\\
			\le &\int_0^c|\phi|'(t)dt\cdot |\nc(K_{\mathcal{L}f}^c)[\omega]|+\left|\int_0^c (\phi'(t) 1_{\{\phi>0\}}(t)-\phi'(t) 1_{\{\phi<0\}(t)})]\nc(K_{\mathcal{L}f}^t)[\omega]dt\right|\\
			\le& \int_0^c|\phi'(t)|dt\sup_{t\in [0,c]}|\nc(K_{\mathcal{L}f}^t)[w]|+\int_0^c |\phi'(t)|\cdot|\nc(K_{\mathcal{L}f}^t)[\omega]|dt\\
			\le&2\int_0^c|\phi'(t)|dt\sup_{t\in [0,c]}|\nc(K_{\mathcal{L}f}^t)[w]|.
		\end{align*}
	\end{proof}
	
	For the construction of the functional intrinsic volumes, we will need to consider slices of the differential cycle. Recall that a current $S$ on an open subset $U$ of $\R^n$ is called \emph{normal} if $S$ and $\partial S$ have finite mass. It is called \emph{locally normal} if the restrictions of $S$ and $\partial S$ to any relatively compact subset of $U$ have finite mass. Let $N_{k,loc}(U)$ denote the space of locally normal $k$-currents on $U$, which becomes a locally convex vector space with respect to the family of semi-norms $S\mapsto M_A(S)+M_A(\partial S)$ for $A\subset U$ compact. Let $S\in N_{k,loc}(U)$ and $h:U\rightarrow\R$ be a smooth map that is proper on $\supp S$. Then there exists a Lebesgue-almost everywhere defined, measurable function $\langle S,h,\cdot\rangle:\R\rightarrow N_{k-1,loc}(U)$ such that for every $\phi\in C_c(\R)$
	\begin{align}
		\label{eq:slicingFormula}
		\int_{\R}\phi(s)\langle S,h,s\rangle (\psi)ds=S( h^*(\phi(s)ds)\wedge\psi)\quad\text{for all } \psi\in \Omega_c^{k-1}(U)
	\end{align}
	compare \cite{FedererGeometricmeasuretheory1969} Theorem 4.3.2. Moreover, this map is unique up to sets of zero Lebesgue measure and given by
	\begin{align*}
		\langle S,h,s\rangle=\partial (S\llcorner1_{\{h<s\}})-\partial S\llcorner 1_{\{h>s\}}=-\partial (S\llcorner1_{\{h>s\}})+\partial S\llcorner1_{\{h<s\}}\quad\text{for almost all } s\in \R.
	\end{align*}
	We will refer to \eqref{eq:slicingFormula} as the slicing formula.\\
	
	For smooth functions, we have the following relation between sublevel sets of $\mathcal{L}f$ and slices of the differential cycle $D(f)$. Let $\pi_2:\R^n\times\R^n\rightarrow\R^n$ denote the projection onto the second factor. We may naturally consider the Legendre transform of $f$ as a function on the cotangent bundle by identifying it with the function $\pi^*_2\mathcal{L}f$. 
	\begin{proposition}
		\label{proposition:connection__differential_cycle_normal_cycle_sublevel_sets}
		Let $f\in \Conv_0^+(\R^n,\R)\cap C^\infty(\R^n)$. Consider the map
		\begin{align*}
			F_f:\R^n\times S^{n-1}&\rightarrow \R^n\times\R^n\\
			(y, v)&\mapsto \left(|d\mathcal{L}f(y)|\cdot v,y\right).
		\end{align*}
		Then $(F_f)_{\#}\nc(K_{\mathcal{L}f}^t)=\langle D(f),\pi_2^*\mathcal{L}f,t\rangle $ for almost all $t>0$.
	\end{proposition}
	\begin{proof}
		As $D(f)$ is closed,
		\begin{align*}
			\langle D(f),\pi_2^*\mathcal{L}f,t\rangle= \partial (D(f)\llcorner 1_{\{\mathcal{L}f< t\}})\quad\text{for almost all }t\in(0,\infty).
		\end{align*}
		But $D(f)\llcorner 1_{\{\mathcal{L}f< t\}}$ is given by integration over the smooth oriented submanifold $\{(x,y):y=df(x),\mathcal{L}f(y)<t\}$, whose boundary is $\{(x,y):y=df(x),\mathcal{L}f(y)=t\}$ (with the orientation induced by the pullback of the $1$-form $\pi_2 ^*d\mathcal{L}f$ to $T^*\R^n$). Stokes theorem implies that $\langle D(f),\pi_2^*\mathcal{L}f,t\rangle$ is given by integration over the smooth submanifold $\{(x,y):y=df(x),\mathcal{L}f(y)=t\}$ for almost all $t\in(0,\infty)$. As $\mathcal{L}f$ is a smooth function, the normal cycle of $K_{\mathcal{L}f}^t$ is given by $\{(y,\frac{d\mathcal{L}f(y)}{|d\mathcal{L}f(y)|}):y\in\R^n\}$ with the orientation induced by the restriction of the $1$-form $d\mathcal{L}f$ to $S\R^n$. Obviously, $F_f$ establishes an orientation preserving diffeomorphism between these two smooth submanifolds. The claim follows.
	\end{proof}

	Finally, let us note the following relation.
	\begin{lemma}
		\label{lemma:relationBetaLegendre}
		Let $f\in \Conv_0^+(\R^n,\R)\cap C^\infty$. Then for every $\omega\in \Omega^{n-1}(T^*\R^n)$ with $\pi(\supp \omega)$ bounded,
		\begin{align*}
			D(f)[\beta\wedge \omega]=D(f)[\pi_2^*d\mathcal{L}f\wedge\omega].
		\end{align*}
	\end{lemma}
	\begin{proof}
		Since $D(f)$ is given by integration over the graph of $df$, it is sufficient to show that the restrictions of $\beta$ and $\pi_2^*d\mathcal{L}f$ coincide. As the equation $f(x)+\mathcal{L}f(y)=\langle y,x\rangle$ holds on the graph of $df$, the pullback of the form $\pi^*df+\pi_2^*d\mathcal{L}f-d\langle y,x\rangle$ to the graph of $df$ vanishes identically. A short calculation shows that $d\langle y,x\rangle =\alpha+\beta$. Noting that $\pi^*df$ and $\alpha$ restrict to the same form on the graph of $df$, we therefore obtain that $\pi^*d\mathcal{L}f$ and $\beta$ also restrict to the same form.
	\end{proof}

	\section{Invariant differential forms and the Hessian measures}
	\label{section:InvDiffFormsAndHessiansMeasures}
	\subsection{The Hessian measures as integrals over the differential cycle}
	\label{section:repHessianMeasures}
	In this section we will show how the Hessian measures can be represented in terms of integrating certain differential forms over the differential cycle. The following result is well known. We include a proof for the convenience of the reader.
	\begin{proposition}
		\label{proposition:rotation-invariant-forms:SO(n)-invariant-forms_linearized_version}
		The forms 
		\begin{align*}
			\omega_s&=\sum\limits_{k=1}^{n}dx_k\wedge dy_k,\\
			\kappa_i&=\frac{1}{i!(n-i)!}\sum\limits_{\sigma\in S_{n}}\sign(\sigma)dx_{\sigma(1)}\dots dx_{\sigma(n-i)}\wedge dy_{\sigma(n-i+1)}\dots dy_{\sigma(n)} \quad 0\le i\le n,
		\end{align*}
		generate the algebra of all $\mathrm{SO}(n)$-invariant forms in $\Lambda^*(\R^n\times\R^n)$.
	\end{proposition}
	\begin{proof}
		Let $\Lambda^*(\R^n\times\R^n)_{\C}^{\SO(n)}$ denote the space of $\mathrm{SO}(n)$-invariant differential forms with complex coefficients and consider the operation of $\mathrm{O}(n)$ on $\Lambda^*(\R^n\times\R^n)_{\C}^{\mathrm{SO}(n)}$. Note that $\mathrm{SO}(n)$ is a normal subgroup of $\mathrm{O}(n)$, so this operation is well defined. In fact, this operation descends to an operation of the group $\mathrm{O}(n)/\mathrm{SO}(n)$, which is isomorphic to the commutative group $\{\pm 1\}$. This implies in particular that the only $\mathrm{O}(n)$-irreducible subspaces are $1$-dimensional. As $\mathrm{O}(n)$ is compact, we can consequently decompose this space into $1$-dimensional, irreducible subspaces. In particular, if we fix any $g_0\in \mathrm{O}(n)\setminus\mathrm{SO}(n)$, then the action of $g_0$ on $\Lambda^*(\R^n\times\R^n)_{\C}^{\SO(n)}$ is diagonizable, so we can decompose $\Lambda^*(\R^n\times\R^n)_{\C}^{\mathrm{SO}(n)}$ into a direct sum of eigenspaces. Since $g_0^2\in\mathrm{SO}(n)$, the corresponding eigenvalues are $\pm 1$. Let $\Lambda^*(\R^n\times\R^n)_{\C}^{\mathrm{SO}(n),\pm}$ denote the corresponding eigenspaces.\\
		By construction, $g\in \mathrm{O}(n)$ operates on $\Lambda^*(\R^n\times\R^n)_{\C}^{\mathrm{SO}(n),-}$ by multiplication with $\det(g)$ and trivially on $\Lambda^*(\R^n\times\R^n)_{\C}^{\mathrm{SO}(n),+}$. For $\tau\in \Lambda^*(\R^n\times\R^n)_{\C}^{\mathrm{SO}(n)}$ define
		\begin{align*}
			A_+\tau:=&\frac{1}{n!}\sum_{\sigma\in S_n} g_\sigma^*\tau,\\
			A_-\tau:=&\frac{1}{n!}\sum_{\sigma\in S_n} \sign(\sigma)g_\sigma^*\tau,
		\end{align*}
		where $g_\sigma\in \mathrm{O}(n)$ denotes the matrix with entries $(g_\sigma)_{jk}=\begin{cases}
			1 & k=\sigma(j)\\
			0 & \text{else} 
		\end{cases}$. It is then easy to see that $A_+$ and $A_-$ are the projections onto $\Lambda^*(\R^n\times\R^n)_{\C}^{\mathrm{SO}(n),\pm}$ respectively.\\
		Let $\tau_+ \in\Lambda^*(\R^n\times\R^n)_{\C}^{\mathrm{SO}(n),+}$ and assume that $\tau_+$ contains a non-trivial contribution of the basis vector $dx_{i_1}\dots dx_{i_k}\wedge dy_{j_1}\dots dy_{j_l}$ for suitable indices $i_1<\dots<i_k,j_1<\dots <j_l$. Since the matrices $g=\mathrm{diag}(1,\dots,1,-1,1,\dots,1)\in \mathrm{O}(n)$ operate by multiplication by $1$ on $\tau_+$, we see that $k=l$ and $i_1=j_1,\dots, i_k=j_l$. It is then easy to check that any such term is mapped by $A_+$ to a power of the symplectic form $\omega_s$.\\
		Similarly, let $\tau_- \in\Lambda^*(\R^n\times\R^n)_{\C}^{\mathrm{SO}(n),-}$ and assume that $\tau_-$ contains a non-trivial contribution of the basis vector $dx_{i_1}\dots dx_{i_k}\wedge dy_{j_1}\dots dy_{j_l}$ for suitable indices $i_1<\dots<i_k,j_1<\dots <j_l$. Since the matrices $g=\mathrm{diag}(1,\dots,1,-1,1,\dots,1)\in \mathrm{O}(n)$ operate by multiplication by $\det(g)=-1$ on $\tau_-$, we see that every index $1,\dots,n$ has to occur precisely once in the set $\{i_1,\dots,i_k,j_1,\dots,j_l\}$. In particular, $k=n-l$. It is now easy to check that $A_-$ maps any form $dx_{i_1}\dots dx_{i_{n-l}}\wedge dy_{j_1}\dots dy_{j_l}$ to a multiple of $\kappa_l$.
	\end{proof}

	For $t\in\R$ let $G_t:T\R^n\rightarrow T\R^n$ be given by $G_t(x,y)=(x,y+tx)$. A short calculation shows
	\begin{lemma}
		\label{lemma:rotation-invariant-valuations:relations-tau}
		\begin{align*}
			G_t^*\kappa_n=\sum_{i=0}^{n}t^{n-i}\kappa_i.
		\end{align*}
	\end{lemma}
	The following result shows that the forms $\kappa_i$ induce the Hessian measures when integrated with respect to the differential cycle.
	\begin{proposition}
		\label{proposition:interpretation-hessian-measures}
		For $f\in\Conv(\R^n,\R) \cap C^2(\R^n)$ and $\phi\in C_c(\R^n)$:
		\begin{align*}
			D(f)[\pi^*\phi\wedge\kappa_i]=\binom{n}{i}\int_{\R^n}\phi(x)[D^2f(x)]_i dx.
		\end{align*}
		Thus $D(f)[1_{\pi^{-1}(B)}\wedge \kappa_i]=\binom{n}{i}\Phi_i(f,B)$ for $f\in\Conv(\R^n,\R)$ and all Borel sets $B\subset\R^n$.
	\end{proposition}
	\begin{proof}
		The second statement follows from the first using the continuity of the Hessian measures $\Phi_i$ (see \cite{ColesantiEtAlHessianvaluations2020} Theorem 7.3) and Corollary \ref{corollary:differential_cycle_vague_continuity}, as well as the representation of the Hessian measures on the dense subspace smooth convex functions, see \cite{ColesantiEtAlHessianvaluations2020} Section 8.5.\\
		Let $h\in\Conv(\R^n,\R)$ be given by $h(x)=\frac{|x|^2}{2}$. Then for $f\in\Conv(\R^n,\R)\cap C^2(\R^n)$ 
		\begin{align*}
			D(f+th)[\pi^*\phi \wedge \kappa_n]=G_{t\#}D(f)[\pi^*\phi \wedge \kappa_n],
		\end{align*}
		for $t\ge0$ by Proposition \ref{proposition_Fu_sum_MA_andC11}. Thus
		\begin{align*}
			D(f+th)[\pi^*\phi \wedge \kappa_n]=&D(f)[\pi^*\phi \wedge G_{t}^*\kappa_n]=\sum_{i=0}^{n}t^{n-i}D(f)[\pi^*\phi \wedge \kappa_i]
		\end{align*}
		by Lemma \ref{lemma:rotation-invariant-valuations:relations-tau}. On the other hand, $\kappa_n=dy_1\wedge\dots\wedge dy_n$, so
		\begin{align*}
			D(f+th)[\pi^*\phi\wedge  \kappa_n]=&\int_{\R^n}\phi(x)\det(D^2f(x)+tId_n)dx
			=\sum_{i=0}^{n}t^{n-i}\binom{n}{i}\int_{\R^n}\phi(x)[D^2f(x)]_idx,
		\end{align*}
		where $Id_n$ denote the identity matrix. Comparing coefficients, we obtain
		\begin{align*}
			D(f)[\pi^*\phi \wedge \kappa_i]=\binom{n}{i}\int_{\R^n}\phi(x)[D^2f(x)]_idx.
		\end{align*}
	\end{proof}

	\subsection{Invariant differential forms}
	
	We define the following $\SO(n)$-invariant differential forms on $T^*\R^n\cong\R^n\times\R^n$
	in degree $1$
	\begin{align*}
		\alpha:=&\sum\limits_{k=1}^{n}y_kdx_k,\\
		\beta:=& \sum\limits_{k=1}^{n}x_kdy_k,\\
		\gamma:=&\sum\limits_{k=1}^{n}x_kdx_k,
	\end{align*}
	and for $0\le i\le n-1$ in degree $n-1$
	\begin{align*}
		\tau_i:=&\frac{1}{i!(n-i-1)!}\sum\limits_{\sigma\in S_n}\sign(\sigma)x_{\sigma(1)}dx_{\sigma(2)}\dots dx_{\sigma(n-i)}\wedge dy_{\sigma(n-i+1)}\dots dy_{\sigma(n)}.
	\end{align*}
	Note that $d\tau_i=(n-i)\kappa_i$ are exactly the invariant forms from Proposition \ref{proposition:rotation-invariant-forms:SO(n)-invariant-forms_linearized_version} and $\alpha$ is the canonical $1$-form on $T^*\R^n$, so that the symplectic form is given by $\omega_s=-d\alpha=\sum_{i=k}^{n}dx_k\wedge dy_k$. Let us compare these forms to the forms $\tilde{\kappa}_i$ from Section \ref{section:IntegralGeo}. If $j:\R^n\times\R^n\rightarrow\R^n\times\R^n$ denotes the map that interchanges the two factors, then
	\begin{align*}
		(j^*\tau_i)|_{S\R^n}=\tilde{\kappa}_i.
	\end{align*}
	From the definition of these forms one easily deduces
	\begin{lemma}
		The forms $\gamma\wedge\tau_i$ and $\beta\wedge\tau_i$ are primitive, that is, their product with $\omega_s$ vanishes.
	\end{lemma}
	\begin{proof}
		As the forms are $\SO(n)$-invariant, it is sufficient to show that the products with $\omega_s$ vanish in $(te_1,0)$ for $t\ge 0$. Note that
		\begin{align*}
			\gamma|_{(te_1,0)}=&tdx_1,\\
			\tau_i=&\frac{1}{i!(n-i-1)!}\sum\limits_{\sigma\in S_n,\sigma(1)=1}\sign(\sigma)tdx_{\sigma(2)}\dots dx_{\sigma(n-i)}\wedge dy_{\sigma(n-i+1)}\dots dy_{\sigma(n)}, 
		\end{align*}
		so $(\omega_s\wedge\gamma\wedge \tau_i)|_{(te_1,0)}$ is proportional to
		\begin{align*}
			\left(\sum_{k=2}^{n}dx_i\wedge dy_i\right)\wedge tdx_1\wedge \sum\limits_{\sigma\in S_n,\sigma(1)=1}\sign(\sigma)tdx_{\sigma(2)}\dots dx_{\sigma(n-i)}\wedge dy_{\sigma(n-i+1)}\dots dy_{\sigma(n)},
		\end{align*}
		were the form $dx_{\sigma(2)}\dots dx_{\sigma(n-i)}\wedge dy_{\sigma(n-i+1)}\dots dy_{\sigma(n)}$ contains one of the forms $dx_k$ or $dy_k$ for each $2\le k\le n$. Thus this term vanishes identically. A similar argument applies to $\beta\wedge\tau$.
	\end{proof}

	Let $r:T\R^n\rightarrow\R$ be the function given by $r(x,y)=|x|$. Then $r^2:T\R^n\rightarrow\R$ is smooth and $d r^2=2\gamma$.
	\begin{proposition}
		\label{proposition:invariant_forms_Relations}
		For $1\le i\le n-1$,
		\begin{align*}
			r^2\kappa_i=\gamma\wedge \tau_i+\beta\wedge\tau_{i-1}.
		\end{align*}
		In addition,
		\begin{align*}
			r^2\kappa_0=&\gamma\wedge\tau_0,\\
			r^2\kappa_n=&\beta\wedge\tau_{n-1}.
		\end{align*}
	\end{proposition}
	\begin{proof}
		Let $X_\beta$ and $X_\gamma$ denote the symplectic vector fields of $\beta$ and $\gamma$, that is, the unique vector fields that satisfy
		\begin{align*}
			&i_{X_\beta}\omega_s=\beta, && i_{X_\gamma}\omega_s=\gamma.
		\end{align*}
		Then 
		\begin{align*}
			&X_\gamma=\sum_{k=1}^{n}-x_k\frac{\partial }{\partial y_k}, 		&&X_\beta=\sum_{k=1}^{n}x_k\frac{\partial }{\partial x_k}.
		\end{align*}
		In particular, $i_{X_{\gamma}}{\beta}=-r^2$, $i_{X_{\beta}}{\gamma}=r^2$, which implies
		\begin{align*}
			i_{X_{\gamma}}({\beta}\wedge\kappa_i)=-r^2\kappa_i-{\beta}\wedge i_{X_{\gamma}}\kappa_i.
		\end{align*}
		On the other hand, $\kappa_i$ is primitive, that is, $\omega_s\wedge\kappa_i=0$, which implies
		\begin{align*}
			i_{X_{\gamma}}({\beta}\wedge\kappa_i)=&i_{X_{\gamma}}(i_{X_{\beta}}\omega_s\wedge\kappa_i)=i_{X_{\gamma}}\left(i_{X_{\beta}}(\omega_s\wedge\kappa_i)-\omega_s\wedge i_{X_{\beta}}\kappa_i\right)\\
			=&-i_{X_{\gamma}}\omega_s\wedge i_{X_{\beta}}\kappa_i-\omega_s\wedge i_{X_{\gamma}}i_{X_{\beta}}\kappa_i=-{\gamma}\wedge i_{X_{\beta}}\kappa_i-\omega_s\wedge i_{X_{\gamma}}i_{X_{\beta}}\kappa_i.
		\end{align*}
		Thus,
		\begin{align*}
			r^2\kappa_i=&-i_{X_{\gamma}}({\beta}\wedge\kappa_i)-{\beta}\wedge i_{X_{\gamma}}\kappa_i\\
			=&{\gamma}\wedge i_{X_{\beta}}\kappa_i+\omega_s\wedge i_{X_{\gamma}}i_{X_{\beta}}\kappa_i-{\beta}\wedge i_{X_{\gamma}}\kappa_i.
		\end{align*}
		It is thus sufficient to establish the formulas $i_{X_{\beta}}\kappa_i=\tau_i$, $i_{X_{\gamma}}i_{X_{\beta}}\kappa_i=0$, and $i_{X_{\gamma}}\kappa_i=-\tau_{i-1}$. Due to the invariance properties of the forms involved, it is furthermore sufficient to establish the equation in $(x,y)=(te_1,0)$ for all $t\ge0$, which is a simple calculation. We leave the details to the reader.
	\end{proof}
	
	Let us remark that this proposition may be seen as version of Equation (11) in \cite{BrandoliniEtAlSerrintypeoverdetermined2008}, see also \cite{ColesantiEtAlHadwigertheoremconvex2020} Lemma 3.10. More precisely, the formula given there can be deduced by combining Proposition \ref{proposition:invariant_forms_Relations} with Lemma \ref{lemma:pullbackTauSurfaceAreaMeasure} below.

	\section{Singular Hessian valuations}
	\label{section:singularHessianValuations}
	The goal of this section is to reinterpret the original construction of the functional intrinsic volumes by Colesanti, Ludwig and Mussnig \cite{ColesantiEtAlHadwigertheoremconvex2020} in terms of the differential cycle. To facilitate a comparison between the two formulations, we will state which of the steps below  have a corresponding result in \cite{ColesantiEtAlHadwigertheoremconvex2020}.\\
	
	In the first subsection, we will prove an inequality for smooth functions in $\Conv_0^+(\R^n,\R)$. In \cite{ColesantiEtAlHadwigertheoremconvex2020}, this corresponds to an inequality for expressions of the form
	\begin{align*}
		\int_{t_1<\mathcal{L}f(d f(x))\le t_2}\zeta(|x|)d\Phi_i(f,x)=\binom{n}{i}^{-1}D(f)[ 1_{\{t_1<\pi_2^*\mathcal{L}f\le t_2\}}\zeta(|x|)\kappa_i],
	\end{align*}
	where $\pi_2:\R^n\times\R^n\rightarrow\R^n$ is the projection onto the second factor, $f\in\Conv_0^+(\R^n,\R)\cap C^\infty(\R^n)$ and $\zeta\in D^n_i$. We will simplify this in two ways: First, it will be sufficient to only derive the estimates for continuous $\zeta\in C_c([0,\infty))$. Secondly, we will replace the indicator $1_{\{t_1<\pi_2^* \mathcal{L}f\le t_2\}}$ by $\phi(\mathcal{L}f(y))$ for $\phi\in C_c^\infty((0,\infty))$. As we will see in the next subsection, we can easily remove this function from the inequality and obtain a form that is valid for all convex functions. 
	
	\subsection{An inequality}
	
	Throughout this section let $1\le i\le n-1$. The goal of this section is the proof of the following inequality, which corresponds to \cite{ColesantiEtAlHadwigertheoremconvex2020} Lemma 3.14. In order to keep the notation concise, we denote the standard coordinates on $\R^n\times\R^n$ by $(x,y)$ and denote the pullback of a function $h$ defined on $\R^n$ along the projections onto the two different factors by $h(x)$ or $h(y)$ respectively.
	\begin{proposition}
		\label{proposition:main_inequality}
		If $f\in\Conv^+_0(\R^n,\R)\cap C^\infty(\R^n)$ and $\zeta\in C_c([0,\infty))$ with $\supp\zeta\subset [0,R]$ for some $R>0$, then the following inequality holds for all $\phi\in C^\infty_c((0,\infty))$:
		\begin{align*}
			&\left|D(f)\left[ \phi(\mathcal{L}f(y))\zeta(|x|) \kappa_i\right]\right|
			\le 2^{i+1}\omega_{n}\binom{n}{i}\left(\sup_{|x|\le R+1}|f(x)|\right)^i\|\zeta\|\int_0^{c(f)}|\phi'(t)|dt,
		\end{align*}
		where $c(f):=(1+2R)\sup_{|x|\le R+1}|f(x)|$.
	\end{proposition}
	
	The result will follow directly by combining Lemma \ref{lemma:estimatePsi} and Lemma \ref{lemma:estimatePsiPrime} below. 
	We will need some preparation.\\
	
	Let us start with the following observation, which corresponds to \cite{ColesantiEtAlHadwigertheoremconvex2020} Proposition 3.12:
	\begin{corollary}
		\label{corollary:differential_gamma_tau_in_terms_of_theta}
		For any $\Psi\in C^1((0,\infty))$
		\begin{align*}
			d(\Psi(|x|)\tau_{i})=&\Psi(|x|)\kappa_i+|x|\Psi'(|x|)\kappa_i-\frac{\Psi'(|x|)}{|x|}\beta\wedge\tau_{i-1}
		\end{align*}
		on $(\R^n\setminus\{0\})\times\R^n$. 
	\end{corollary}
	\begin{proof}
		Using $d\tau_i=(n-i)\kappa_i$, this follows from Proposition \ref{proposition:invariant_forms_Relations}:
		\begin{align*}
			d(\Psi(|x|)\tau_{i})=&(n-i)\Psi(|x|)d\tau_{i}+\frac{\Psi'(|x|)}{|x|}\gamma\wedge\tau_{i}\\
			=&(n-i)\Psi(|x|)\kappa_i+\frac{\Psi'(|x|)}{|x|}\left[|x|^2\kappa_i-\beta\wedge\tau_{i-1}\right]
		\end{align*}
		on $(\R^n\setminus\{0\})\times\R^n$.
	\end{proof}

	For the rest of this subsection, $\zeta\in C_c([0,\infty))$ will be a fixed function with $\supp\zeta\subset [0,R]$. We will associate the two continuous functions $\eta,\rho\in C_c([0,\infty))$ to $\zeta$ defined for $t\ge 0$ by
	\begin{align*}
		&\eta(t):=\int_t^\infty r^{n-i-1}\zeta(r)dr &&\rho(t):=t^{n-i}\zeta(t)+(n-i)\eta(t).
	\end{align*}
	Note that $\supp\zeta\subset [0,R]$ implies that $\supp\eta,\supp\rho\subset [0,R]$.\\
	
	If we choose (as in the proof of \cite{ColesantiEtAlHadwigertheoremconvex2020} Proposition 3.13)
	\begin{align*}
		\Psi(r):=-\frac{1}{r^{n-i}}\int_r^\infty s^{n-i-1}\zeta(s)ds,
	\end{align*}
	then $\Psi$ is a $C^1$-function on $(0,\infty)$ with
	\begin{align*}
		\Psi'(r)=\frac{n-i}{r^{n-i+1}}\int_r^\infty s^{n-i-1}\zeta(s)ds+\frac{\zeta(r)}{r},
	\end{align*}
	so
	\begin{align*}
		(n-i)\Psi(r)+r\Psi'(r)=\zeta(r)
	\end{align*}
	and
	\begin{align*}
		r^{n-i+1}\Psi'(r)=(n-i)\eta(r)+r^{n-i}\zeta(r)=\rho(r).
	\end{align*}
	In particular, we may apply Corollary \ref{corollary:differential_gamma_tau_in_terms_of_theta} to obtain
	\begin{align}
		\notag
		\zeta(|x|)\kappa_i=&\left[(n-i)\Psi(|x|)+|x|\Psi'(|x|)\right]\kappa_i\\
		\label{equation:relation_zeta_Theta_vs_gamma_tau}
		=&d(\Psi(|x|)\tau_{i})+\frac{\Psi'(|x|)}{|x|}\beta\wedge\tau_{i-1}\quad\text{on }(\R^n\setminus\{0\})\times\R^n.
	\end{align}

	Let us now turn to the proof of Proposition \ref{proposition:main_inequality}. First note that the restrictions of $\beta$ and $\pi_2^*d\mathcal{L}f$ coincide on the graph of $df$ for any $f\in\Conv_0^+(\R^n ,\R)\cap C^\infty(\R^n )$ by Lemma \ref{lemma:relationBetaLegendre}. Equation \eqref{equation:relation_zeta_Theta_vs_gamma_tau} thus implies for $\phi\in C^\infty_c((0,\infty))$
	\begin{align*}
		&D(f)[\phi(\mathcal{L}f(y))\zeta(|x|)\kappa_i)]\\
		=&D(f)\left[\phi(\mathcal{L}f(y))\left(d(\Psi(|x|)\tau_{i})+\frac{\Psi'(|x|)}{|x|}\beta\wedge\tau_{i-1}\right)\right]\\
		=&-D(f)[\phi'(\mathcal{L}f(y))d\pi_2^*\mathcal{L}f\wedge \Psi(|x|)\tau_{i}]+D(f)\left[\phi(\mathcal{L}f(y))\frac{\Psi'(|x|)}{|x|}\beta\wedge\tau_{i-1}\right]\\
		=&-\int_0^\infty \phi'(t)\langle D(f),\pi_2^*\mathcal{L}f,t\rangle [\Psi(|x|)\tau_{i}]dt+\int_0^\infty \phi(t)\langle D(f),\pi_2^*\mathcal{L}f,t\rangle \left[\frac{\Psi'(|x|)}{|x|}\tau_{i-1}\right]dt,
	\end{align*}	
	where we have used the slicing formula for each of the two terms and that $D(f)$ is closed.\\
	
	We will need the following relation between the forms considered in the previous section involving the map $F_f:\R^n\times S^{n-1}\rightarrow \R^n\times\R^n$ considered in Proposition \ref{proposition:connection__differential_cycle_normal_cycle_sublevel_sets}. 
	\begin{lemma}
		\label{lemma:pullbackTauSurfaceAreaMeasure}
		\begin{align*}
			F_\mathcal{L}f\tau_i=|d\mathcal{L}f(y)|^{n-i}(j^*\tau_i)|_{S\R^n}=|d\mathcal{L}f(y)|^{n-i}\tilde{\kappa}_i
		\end{align*}
	\end{lemma}
	\begin{proof}
		Consider the map 
		\begin{align*}
			G:\R\times \R^n\times \R^n&\rightarrow\R^n\times\R^n\\
			(t,x,y)&\mapsto (tx,y).
		\end{align*}
		From the definition of $\tau_i$, one directly obtains $G^*\tau_i=t^{n-i}\tau_i$. If we consider 
		\begin{align*}
			F:\R^n\times S^{n-1}\rightarrow&\R\times\R^n\times \R^n\\
			(y,v)\mapsto& (|d\mathcal{L}f(y)|,v,y),
		\end{align*}
		then $F_f=G\circ F$, which implies the claim.
	\end{proof}
	
	\begin{lemma}
		\label{lemma:estimatePsi}
		For $f\in\Conv_0^+(\R^n,\R)\cap C^\infty(\R^n)$, $\phi\in C^\infty_c((0,\infty))$:
		\begin{align*}
			&\left|\int_0^\infty \phi'(t)\langle D(f),\pi_2^*\mathcal{L}f,t\rangle [\Psi(|x|)\tau_{i}]dt\right|\\
			\le&  2^i(n-i)\omega_{n-i}\mu_i(B_1(0))\left(\sup_{|x|\le R+1}|f(x)|\right)^i\|\eta\|_\infty\int_0^{c(f)}|\phi'(t)|dt,
		\end{align*}
		where $c(f)$ is the constant defined in Proposition \ref{proposition:main_inequality}.
	\end{lemma}
	\begin{proof}
		Using Proposition \ref{proposition:connection__differential_cycle_normal_cycle_sublevel_sets} and Lemma \ref{lemma:pullbackTauSurfaceAreaMeasure}, we obtain
		\begin{align*}
			\int_0^\infty \phi'(t)\langle D(f),\pi_2^*\mathcal{L}f,t\rangle [\Psi(|x|)\tau_{i}]dt=&\int_0^\infty\phi'(t)\nc(K_{\mathcal{L}f}^t)\left[|d\mathcal{L}f(y)|^{n-i}\Psi(|d\mathcal{L}f(y)|)(j^*\tau_i)|_{S\R^n}\right]dt\\
			=&-\int_0^\infty\phi'(t)\nc(K_{\mathcal{L}f}^t)\left[\eta(|d\mathcal{L}f(y)|)\tilde{\kappa}_i\right]dt.
		\end{align*}	
		Note that $\eta$ is supported on $[0,R]$. If $|d\mathcal{L}f(y)|\le R$, then $y= df(x_0)$ for $x_0=d\mathcal{L}f(y)\in\R^n$ with $ |x_0|\le R$, so Lemma \ref{lemma:bounds_subgradient} implies $|\mathcal{L}f(y)|\le c(f)$ and $|y|\le 2\sup_{|x|\le R+1}|f(x)|$. 	As $\tilde{\kappa}_i$ induces a positive curvature measure, this implies
		\begin{align*}
			\left|\int_0^\infty \phi'(t)\langle D(f),\pi_2^*\mathcal{L}f,t\rangle [\Psi(|x|)\tau_{i}]dt\right|\le& \|\eta\|_\infty\int_0^{c(f)}|\phi'(t)| \nc(K_{\mathcal{L}f}^t)[\tilde{\kappa}_i]dt\\
			=&(n-i)\omega_{n-i}\|\eta\|_\infty\int_0^{c(f)}|\phi'(t)| \mu_i(K_{\mathcal{L}f}^t)dt\\
			\le&(n-i)\omega_{n-i}\|\eta\|_\infty\int_0^{c(f)}|\phi'(t)|dt\cdot  \mu_i(K_{\mathcal{L}f}^{c(f)}).
		\end{align*}
		Here we have used that $\mu_i$ is monotone. As $y\in K_{\mathcal{L}f}^t$ implies $|y|\le 2\sup_{|x|\le R+1}|f(x)|$ by Lemma \ref{lemma:bounds_subgradient}, $K_{\mathcal{L}f}^t\subset 2\sup_{|x|\le R+1}|f(x)|\cdot B_1(0)$. The monotonicity of the $i$th intrinsic volume thus implies
		\begin{align*}
			&\left|\int_0^\infty \phi'(t)\langle D(f),\pi_2^*\mathcal{L}f,t\rangle [\Psi(|x|)\tau_{i}]dt\right|\\
			\le& 2^i(n-i)\omega_{n-i}\mu_i(B_1(0))\left(\sup_{|x|\le R+1}|f(x)|\right)^i\|\eta\|_\infty\int_0^{c(f)}|\phi'(t)|dt.
		\end{align*}
	\end{proof}
	
	\begin{lemma}
		\label{lemma:estimatePsiPrime}
		For $f\in\Conv_0^+(\R^n,\R)\cap C^\infty(\R^n)$, $\phi\in C^\infty_c((0,\infty))$:
		\begin{align*}
			&\left|\int_0^\infty \phi(t)\langle D(f),\pi_2^*\mathcal{L}f,t\rangle [\Psi'(|x|)\tau_{i-1}]dt\right|\\
			\le&2^{i+1}\omega_{n-i}\mu_i(B_1(0))\left(\sup_{|x|\le R+1}|f(x)|\right)^i\|\rho\|_\infty\int_0^{c(f)}|\phi'(t)|dt,
		\end{align*}
		where $c(f)$ is the constant defined in Proposition \ref{proposition:main_inequality}.
	\end{lemma}
	\begin{proof}
		We use the same reasoning as in the proof of Lemma \ref{lemma:estimatePsi}. First we use the relation $(F_f)_{\#}\nc(K_{\mathcal{L}f}^t)=\langle D(f),\pi_2^*\mathcal{L}f,t\rangle $ for almost all $t>0$ to obtain
		\begin{align*}
			\int_0^\infty \phi(t)\langle D(f),\pi_2^*\mathcal{L}f,t\rangle \left[\frac{\Psi'(|x|)}{|x|}\tau_{i-1}\right]dt=&\int_0^\infty \phi(t)\nc(K_{\mathcal{L}f}^t) \left[\frac{\rho(|d\mathcal{L}f(y)|)}{|d\mathcal{L}f(y)|}\tilde{\kappa}_{i-1}\right]dt\\
			=&\frac{1}{n-i}\int_0^\infty \phi(t)\nc(K_{\mathcal{L}f}^t) \left[\frac{\rho(|d\mathcal{L}f(y)|)}{|d\mathcal{L}f(y)|}i_TD\tilde{\kappa}_{i}\right]dt,
		\end{align*}
		where we have used Equation \eqref{eq:RuminDiffTildeKappa} in the last step. Note that $\eta$ is supported on $[0,R]$ by assumption and as in the proof of Lemma \ref{lemma:estimatePsi}, $y\in \partial K_{\mathcal{L}f}(t)$ and $|d\mathcal{L}f(y)|\le R$ imply $|\mathcal{L}f(y)|\le c(f)$. As $i_TD\tilde{\kappa}_i$ defines a non-negative measure, we thus obtain
		\begin{align*}
			\left|\int_0^\infty \phi(t)\langle D(f),\pi_2^*\mathcal{L}f,t\rangle \left[\frac{\Psi'(|x|)}{|x|}\tau_{i-1}\right]dt\right|\le&\frac{1}{n-i}\|\rho\|_\infty\int_0^{c(f)} |\phi(t)|\nc(K_{\mathcal{L}f}^t) \left[\frac{1}{|d\mathcal{L}f(y)|}i_TD\tilde{\kappa}_i\right]dt.
		\end{align*}
		Corollary \ref{corollary:inequality_partial_integration} and the monotonicity of $\mu_i$ imply
		\begin{align*}
			\int_0^{c(f)} |\phi(t)|\nc(K_{\mathcal{L}f}^t) \left[\frac{1}{|d\mathcal{L}f(y)|}i_TD\tilde{\kappa}_i\right]dt\le& 2\int_0^{c(f)}|\phi'(t)|dt\cdot \nc(K_{\mathcal{L}f}^{c(f)})[\tilde{\kappa}_i]\\
			=&2(n-i)\omega_{n-i}\int_0^{c(f)}|\phi'(t)|dt\cdot \mu_i(K_{\mathcal{L}f}^{c(f)}).
		\end{align*}
		By the same argument as in the proof of Lemma \ref{lemma:estimatePsi}, $K_{\mathcal{L}f}^{c(f)}\subset2\sup_{|x|\le  R+1}|f(x)|\cdot B_1(0)$, so the monotonicity of the $i$th intrinsic volume implies
		\begin{align*}
			&\left|\int_0^\infty \phi(t)\langle D(f),\pi_2^*\mathcal{L}f,t\rangle \left[\frac{\Psi'(|x|)}{|x|}\tau_{i-1}\right]dt\right|\\
			\le&2^{i+1}\omega_{n-i}\mu_i(B_1(0))\left(\sup_{|x|\le R+1}|f(x)|\right)^i\|\rho\|_\infty\int_0^{c(f)}|\phi'(t)|dt.
		\end{align*}
	\end{proof}
	\begin{proof}[Proof of Proposition \ref{proposition:main_inequality}]
		This follows by combining Lemma \ref{lemma:estimatePsi} and Lemma \ref{lemma:estimatePsiPrime} using that  $\mu_i(B_1(0))=\binom{n}{i}\frac{\omega_n}{\omega_{n-i}}$ and $\|\zeta\|=(n-i)\|\eta\|_\infty+\|\rho\|_\infty$.
	\end{proof}
	
	\subsection{Functional intrinsic volumes and principal value integrals}
	\label{Section:functionalIntrinsicVolumes}
	Next, we are going to bring the inequality in Proposition \ref{proposition:main_inequality} into a form that does not depend on $\phi\in C^\infty_c((0,\infty))$ anymore. This will allow us to define the functional intrinsic volumes by continuously extending the map $\zeta\mapsto V^*_{i,\zeta}$ from $C_c([0,\infty))$ to $D^n_i$.
	\begin{proposition}
		\label{prop:inequFuncIntVol}
		Let $1\le i\le n-1$. For every $\zeta\in C_c([0,\infty))$ with $\supp\zeta\subset [0,R]$ for $R>0$ and every $f\in \Conv(\R^n,\R)$
		\begin{align*}
			|D(f)[\zeta(|x|)\kappa_i]|\le 2^{3i+1}\omega_n\binom{n}{i}\left(\sup_{|x|\le R+1}|f(x)|\right)^i\|\zeta\|.
		\end{align*}
	\end{proposition}
	\begin{proof}
		For $f\in\Conv_0^+(\R^n,\R)\cap C^\infty(\R^n)$, Proposition \ref{proposition:main_inequality} implies for $\phi\in C^\infty_c(0,\infty)$
		\begin{align*}
			\left|D(f)\left[ \phi(\mathcal{L}f(y))\zeta(|x|) \kappa_i\right]\right|
			\le 2^{i+1}\omega_n\binom{n}{i}\left(\sup_{|x|\le R+1}|f(x)|\right)^i\int_0^{c(f)}|\phi'(t)|dt \cdot \|\zeta\|.
		\end{align*}
		Let $\psi_1\in C^\infty(0,\infty)$ be a function with $\psi_1=0$ on $[0,\frac{1}{2}]$, $\psi_1=1$ on $[1,\infty)$ and let $\psi_2\in C^\infty(0,\infty)$ be a function with $\psi_2=1$ on $(0,1]$, $\psi_2=0$ on $[2,\infty)$. We may assume that $\psi_1$ is non-decreasing. For $\delta\in (0,1)$ we define $\psi_\delta\in C^\infty_c(0,\infty)$ by
		\begin{align*}
			\psi_\delta(t):=\psi_1\left(\frac{t}{\delta}\right)\psi_2(\delta t).
		\end{align*}
		Then $\psi_\delta\equiv 1$ on $[\delta,\frac{1}{\delta}]$. Moreover,
		\begin{align*}
			\lim\limits_{\delta\rightarrow0}\int_0^{c(f)}|\psi'_\delta(t)|ds=\lim\limits_{\delta\rightarrow 0}\int_0^{c(f)}\frac{1}{\delta}\left|\psi_1'(\frac{t}{\delta})\right|dt=\lim\limits_{\delta\rightarrow 0}\int_0^{\frac{c(f)}{\delta}}\psi_1'(s)ds=\psi_1(1)-\psi_1(0)=1,
		\end{align*}
		because $\psi_1$ is non-decreasing and $\psi_1'$ is supported on $[0,1]$. Next, note that because $f\in \Conv_0^+(\R^n,\R)\cap C^\infty(\R^n)$, $\mathcal{L}f$ is a smooth function and because $D(f)$ is given by integration over the graph of $df$, we obtain $\mathcal{L}f(y)=\langle df(x),x\rangle-f(x)$ for all $(x,y)\in\supp D(f)$. Thus
		\begin{align*}
			D(f)[\psi_\delta(\mathcal{L}f(y))\zeta(|x|)\kappa_i]=&D(f)[\psi_{\delta}(\langle df(x),x\rangle-f(x))\zeta(|x|)\kappa_i]\\
			=&\int_{\R^n}[\psi_{\delta}(\langle df(x),x\rangle-f(x))\zeta(|x|)]d\Phi_i(f,x).
		\end{align*}
		Moreover, $\psi_{\delta}(\langle df(x),x\rangle-f(x))$ converges to $1$ for $\delta\rightarrow0$ for all $x\in\R^n$ with $\langle df(x),x\rangle-f(x)\ne0$, that is, for all $x\ne 0$. As $f$ is smooth, the measure $\Phi_i(f)$ is absolutely continuous with respect to the Lebesgue measure, so $x\mapsto \psi_\delta(\langle df(x),x\rangle-f(x))\zeta(|x|)$ converges $\Phi_i(f)$-almost everywhere to $\zeta(|x|)$. As $\zeta(|x|)$ is continuous and compactly supported, it is integrable with respect to $\Phi_i(f)$, so dominated convergence implies
		\begin{align*}
			|D(f)[\zeta(|x|)\kappa_i]|=&\left|\int_{\R^n}\zeta(|x|)d\Phi_i(f,x)\right|
			=\lim\limits_{\delta\rightarrow0}\left|\int_{\R^n}[\psi_{\delta}(\langle df(x),x\rangle-f(x))\zeta(|x|)]d\Phi_i(f,x)\right|\\
			=&\lim\limits_{\delta\rightarrow0}|D(f)[\psi_\delta(\mathcal{L}f(y)\zeta(|x|)\kappa_i)]|\\
			\le&\limsup_{\delta\rightarrow0}2^{i+1}\omega_n\binom{n}{i}\left(\sup_{|x|\le R+1}|f(x)|\right)^i\cdot \|\zeta\| \int_0^{c(f)}|\psi_\delta'(t)|dt\\
			=&2^{i+1}\omega_n\binom{n}{i}\left(\sup_{|x|\le R+1}|f(x)|\right)^i\cdot \|\zeta\|.
		\end{align*}
		Next, let $f\in\Conv(\R^n,\R)$ be an arbitrary smooth function. Then $\tilde{f}:=f-f(0)-\langle d f(0),\cdot\rangle$ belongs to $\Conv_0(\R^n,\R)$. By Proposition \ref{proposition:bound_lipschitz_constant}, $|df(0)|\le \frac{2}{R+1}\sup_{|x|\le R+1}|f(x)|$. As $\tilde{f}+\lambda|\cdot|^2\in\Conv_0^+(\R^n,\R)$ is smooth, we obtain
		\begin{align*}
			|D(f)[\zeta(|x|)\kappa_i]|=&\lim\limits_{\lambda\rightarrow0}|D(\tilde{f}+\lambda|\cdot|^2)[\zeta(|x|)\kappa_i]|\\
			\le&2^{i+1}\omega_n\binom{n}{i}\limsup_{\lambda\rightarrow0}\left(\sup_{|x|\le R+1}|f(x)-f(0)-\langle d f(0),x\rangle +\lambda|x|^2|\right)^i\cdot \|\zeta\|\\
			\le&2^{i+1}\omega_n\binom{n}{i}\limsup_{\lambda\rightarrow0}\left(4\sup_{|x|\le R+1}|f(x)|+\lambda(R+1)^2)\right)^i\cdot \|\zeta\|.
		\end{align*}
		Here we have used that $\zeta(|x|)\kappa_i$ induces a continuous dually epi-translation invariant valuation in the first step. Thus 
		\begin{align*}
			|D(f)[\zeta(|x|)\kappa_i]|\le 2^{3i+1}\omega_n\binom{n}{i}\left(\sup_{|x|\le R+1}|f(x)|\right)^i\|\zeta\|\quad\text{for all } f\in\Conv(\R^n,\R)\cap C^\infty(\R^n).
		\end{align*}
		As both sides depend continuously on $f$ by Corollary \ref{corollary:differential_cycle_vague_continuity}, we obtain the desired inequality.
	\end{proof}
	Recall that $D(f)[1_{\pi^{-1}(B)}\wedge \kappa_i]=\binom{n}{i}\Phi_i(f,B)$ for $f\in\Conv(\R^n,\R)$ and Borel sets $B\subset\R^n$ by Proposition \ref{proposition:interpretation-hessian-measures}.
	\begin{corollary}
		\label{cor:Extension_Vk}
		For $1\le i\le n-1$ there exists a unique continuous map $V^*_{i}:D^{n}_{i,R}\rightarrow\VConv_{i}(\R^n)^{\SO(n)}\cap \VConv_{B_R(0)}(\R^n)$ such that
		\begin{align*}
			V^*_{i}(\zeta)=\int_{\R^n}\zeta(|x|)d\Phi_i(\cdot,x)
		\end{align*}
		for all $\zeta\in C_c([0,\infty))\cap D^{n}_{i,R}$. Moreover,
		\begin{align}
			\label{eq:inequNormVi}
			|V^*_i(\zeta)[f]|\le 2^{3i+1}\omega_n\left(\sup_{|x|\le R+1}|f(x)|\right)^i\|\zeta\|
		\end{align}
		for all $\zeta\in D^n_{i,R}$, $f\in\Conv(\R^n,\R)$.
	\end{corollary}
	\begin{proof}
		Using Proposition \ref{proposition:supportValuation}, it is easy to see that $\supp\zeta\subset[0,R]$ implies that the valuation $f\mapsto\int_{\R^n}\zeta(|x|)d\Phi_i(f,x)$ is supported on $B_R(0)$ for $\zeta\in C_c([0,\infty))$. In particular, the map
		\begin{align*}
			V^*_{i}:C_c([0,\infty))\cap D^{n}_{i,R}&\rightarrow\VConv_{i}(\R^n)^{\SO(n)}\cap \VConv_{B_R(0)}(\R^n)\\
			\zeta&\mapsto \int_{\R^n}\zeta(|x|)d\Phi_i(\cdot,x)
		\end{align*}
		is well defined. Proposition \ref{prop:inequFuncIntVol} shows that it satisfies
		\begin{align*}
			\|V^*_{i}(\zeta)\|_{B_R(0),1}=\sup\left\{|V^*_{i}(\zeta)[f]|:f\in\Conv(\R^n,\R), \ \sup_{|x|\le R+1}|f(x)|\le 1\right\}\le C_{n,i}\|\zeta\|
		\end{align*}
		for $C_{n,i}:=2^{3i+1}\omega_n$ and is thus continuous if we equip $C_c([0,\infty))\cap D^{n}_{i,R}$ with the norm $\|\cdot\|$. By Theorem \ref{thm:BanachStructureCompact}, $\VConv_{i}(\R^n)\cap \VConv_{B_R(0)}(\R^n)$ is a Banach space with respect to the norm $\|\cdot\|_{B_R(0),1}$. Obviously, $\VConv_{i}(\R^n)^{\SO(n)}\cap \VConv_{B_R(0)}(\R^n)$ is a closed subspace and thus a Banach space with respect to this norm.  Thus this map extends uniquely by continuity to the closure of $C_c([0,\infty))\cap D^{n}_{i,R}\subset D^{n}_{i,R}$, which equals $D^{n}_{i,R}$ by Lemma \ref{lemma:density_cont_Dnk}.\\
		
		To see that inequality \eqref{eq:inequNormVi} holds, note that both sides are continuous in $\zeta\in D^n_{i,R}$ with respect to $\|\cdot\|$. As the inequality holds on the dense subspace $C_c([0,\infty))\cap D^n_{i,R}$ by Proposition \ref{prop:inequFuncIntVol}, the claim follows. 
	\end{proof}
	Note that the inclusion $D^n_{i,R'}\rightarrow D^n_{i,R}$ is continuous for $R'<R$. In particular, the valuation $V^*_i(\zeta)$ only depends on $\zeta$ and not the choice of $R>0$ with $\zeta\in D^n_{i,R}$. The valuation $V^*_i(\zeta)$ is thus well defined and independent of this choice. \\
	
	Note that we can combine Lemma \ref{lemma:density_cont_Dnk} and Corollary \ref{cor:Extension_Vk} to obtain the following approximation result, which strengthens \cite{ColesantiEtAlHadwigertheoremconvex2020} Lemma 3.20.
	\begin{corollary}
		\label{cor:approxVkq}
		For $\zeta\in D^{n}_i$ and $r>0$ define $\zeta^r\in C_c([0,\infty))$ by
		\begin{align*}
			\zeta^r(t):=\begin{cases}
				\zeta(t) & \text{for }t> r,\\
				\zeta(r) & \text{for }0\le t\le r.
			\end{cases}
		\end{align*}
		Then $V^*_{i}(\zeta)=\lim\limits_{r\rightarrow0}V^*_i(\zeta^r)$ in $\VConv_i(\R^n)$.
	\end{corollary}
	
We are now going to establish the principal value representation for the functionals $V^*_i(\zeta)$. Let us prove the following slightly refined version of Theorem \ref{maintheorem:PrincipalValue}:
\begin{theorem}
	\label{theorem:PrincipalValueRefined}
	Let $1\le i\le n-1$. For $\zeta\in D^n_i$ and $f\in\Conv(\R^n,\R)$,
	\begin{align*}
		V^*_{i}(\zeta)[f]=\lim\limits_{\epsilon\rightarrow0}\int_{\R^n\setminus B_\epsilon(0)}\zeta(|x|)d\Phi_i(f,x).
	\end{align*}
	Moreover, the convergence is uniform on $\{f\in\Conv(\R^n,\R): \sup_{|x|\le r} |f(x)|\le M\}$ for every $r>1$, $M>0$. In particular, it is uniform on compact subsets of $\Conv(\R^n,\R)$.
\end{theorem}
\begin{proof}
	Consider the function $\zeta^\epsilon$ defined in Lemma \ref{lemma:density_cont_Dnk}. As $\zeta(t)=\zeta^\epsilon(t)$ for $t\ge\epsilon$, we obtain
	\begin{align*}
		&\left|\int_{\R^n}\zeta^\epsilon(|x|)d\Phi_i(f,x)-\int_{\R^n\setminus B_\epsilon(0)}\zeta(|x|)d\Phi_i(f,x)\right|=\left|\int_{B_\epsilon(0)}\zeta^\epsilon(|x|)d\Phi_i(f,x)\right|\\
		\le&|\zeta(\epsilon)|\int_{B_\epsilon(0)}d\Phi_i(f,x)\le |\zeta(\epsilon)|\int_{\R^n}\phi_\epsilon(|x|)d\Phi_i(f,x),
	\end{align*}
	where $\phi_\epsilon:=\phi\left(\frac{|x|}{\epsilon}\right)$ for $\phi\in C_c([0,\infty),[0,1])$ with $\phi\equiv 1$ on $[0,1]$, $\supp\phi\subset [0,2]$. Note that by Proposition \ref{prop:inequFuncIntVol}
	\begin{align*}
		\int_{\R^n}\phi_\epsilon(|x|)d\Phi_i(f,x)\le 2^{3i+1}\omega_n\left(\sup_{|x|\le 2\epsilon+1}|f(x)|\right)^i\|\phi_\epsilon\|.
	\end{align*}
	By a change of variables, we obtain
	\begin{align*}
		\int_t^\infty \phi_\epsilon(s)s^{n-i-1}ds=\int_t^\infty \phi\left(\frac{s}{\epsilon}\right)s^{n-i-1}ds=\epsilon^{n-i}\int_{\frac{t}{\epsilon}}^\infty\phi(s)s^{n-i-1}ds.
	\end{align*}
	Moreover, 
	\begin{align*}
		\sup_{t>0}|t^{n-i}\phi_\epsilon(t)|=\sup_{t>0}\left|t^{n-i}\phi\left(\frac{t}{\epsilon}\right)\right|=\epsilon^{n-i} \sup_{t>0} |t^{n-i}\phi(t)|.
	\end{align*}
	Thus
	\begin{align*}
		\|\phi_\epsilon\|\le& \sup_{t>0}\left|t^{n-i}\phi_\epsilon(t)\right|+2(n-i)\sup_{t>0}\left|\int_t^\infty\phi_\epsilon(s)s^{n-i-1}ds\right|\\
		=& \epsilon^{n-i}\left(\sup_{t>0}\left|t^{n-i}\phi(t)\right|+2(n-i)\sup_{t>0}\left|\int_t^\infty\phi(s)s^{n-i-1}ds\right|\right)\le \epsilon^{n-i}4\|\phi\|.
	\end{align*}
	Combining this estimate with the inequality in Corollary \ref{cor:Extension_Vk} and using that $\zeta-\zeta^\epsilon\in D^n_{i,\epsilon}$,
	\begin{align*}
		&\left|V^*_{i}(\zeta)[f]-\int_{\R^n\setminus B_\epsilon(0)}\zeta(|x|)d\Phi_i(f,x)\right|\\
		\le&|V^*_{i}(\zeta)[f]-V^*_{i}(\zeta^\epsilon)[f]| +\left|\int_{\R^n}\zeta^\epsilon(|x|)d\Phi_i(f,x)-\int_{\R^n\setminus B_\epsilon(0)}\zeta(|x|)d\Phi_i(f,x)\right|\\
		\le &2^{3i+1}\omega_n\left(\sup_{|x|\le 2\epsilon+1}|f(x)|\right)^i\left(\|\zeta-\zeta^\epsilon\|+4\|\phi\|\epsilon^{n-i}|\zeta(\epsilon)|\right).
	\end{align*}
	As $\zeta^\epsilon\rightarrow\zeta$ with respect to $\|\cdot\|$ by Lemma \ref{lemma:density_cont_Dnk} and $\lim\limits_{\epsilon\rightarrow0}\epsilon^{n-i}\zeta(\epsilon)=0$ as $\zeta\in D^n_i$, we see that 
	\begin{align*}
		\lim\limits_{\epsilon\rightarrow0}\left|V^*_{i}(\zeta)[f]-\int_{\R^n\setminus B_\epsilon(0)}\zeta(|x|)d\Phi_i(f,x)\right|=0
	\end{align*}
	uniformly on $\{f\in\Conv(\R^n,\R): \sup_{|x|\le 1+r'}|f(x)|<M\}$ for every $r'>0$, $M>0$. As any compact subset in $\Conv(\R^n,\R)$ is contained in such a subset, the convergence is uniform on compact subsets.
\end{proof}

Although we have associated some valuation $V^*_i(\zeta)\in\VConv_i(\R^n)^{\SO(n)}$ to any $\zeta\in D^n_i$, we still have to show that this functional satisfies the defining property of the functional intrinsic volume $V^*_{i,\zeta}$, that is, that it has the desired representation on convex functions of class $C^2$. We will use Theorem \ref{theorem:PrincipalValueRefined} give a precise description under which circumstances $V^*_i(\zeta)$ admits an integral representation. We start with the following observation. 
\begin{corollary}
	Let $0\le i\le n-1$, $f\in \Conv(\R^n,\R)$. Then $\Phi_i(f)$ is non-atomic. In particular, $\Phi_i(f)[\{0\}]=0$.
\end{corollary}
\begin{proof}
	For $i=0$, $\Phi_0(f)$ is just the Lebesgue measure, which is non-atomic. Let us thus assume that $1\le i\le n-1$.\\
	Note that $\Phi_i$ is equivariant with respect to translations, that is, $\Phi_i(f(\cdot+x))[A]=\Phi_i(f)[A+x]$ for any $f\in\Conv(\R^n,\R)$, $A\subset \R^n$ Borel set. It is thus sufficient to show that $\Phi_i(f)[\{0\}]=0$ for any $f\in\Conv(\R^n,\R)$. As in the proof of Theorem \ref{theorem:PrincipalValueRefined}, we choose $\phi\in C_c([0,\infty))$ with $\supp\phi\subset [0,2]$, $\phi=1$ on $[0,1]$, and consider the function $\phi_\epsilon(t):=\phi(\frac{t}{\epsilon})$. Dominated convergence implies
	\begin{align*}
		\Phi_i(f)[\{0\}]=\lim\limits_{\epsilon\rightarrow0}\int_{\R^n}\phi_\epsilon(|x|) d\Phi_i(f,x)\le \limsup_{\epsilon\rightarrow0}2^{3i+1}\left(\sup_{|x|\le 3}|f(x)|\right)^i\|\phi_\epsilon\|.
	\end{align*}
	As in the proof of Theorem \ref{theorem:PrincipalValueRefined}, $\|\phi\|_\epsilon\le \epsilon^{n-i}4\|\phi\|$, so this limit is $0$.
\end{proof}

The following result shows that $V^*_i(\zeta)$ coincides with the functional intrinsic volume $V^*_{i,\zeta}$.
\begin{corollary}
	\label{corollary:RepresentationFunctIntrinsicVolume}
	Let $1\le i\le n-1$, $\zeta\in D^n_i$, $f\in\Conv(\R^n,\R)$. If $x\mapsto \zeta(|x|)$ is integrable with respect to $\Phi_i(f)$, then
	\begin{align*}
		V^*_i(\zeta)[f]=\int_{\R^n}\zeta(|x|)d\Phi_i(f,x).
	\end{align*}
	This is in particular the case, if $f$ is of class $C^2$ on a neighborhood of $0$.
\end{corollary}
\begin{proof}
	If $x\mapsto \zeta(|x|)$ is integrable with respect to $\Phi_i(f)$, then Theorem \ref{theorem:PrincipalValueRefined} and dominate convergence imply
	\begin{align*}
		V^*_i(\zeta)[f]=\lim\limits_{\epsilon\rightarrow0}\int_{\R^n\setminus B_\epsilon(0)}\zeta(|x|)d\Phi_i(f,x)=\int_{\R^n}\zeta(|x|)d\Phi_i(f,x),
	\end{align*}
	which shows the first claim. If $f$ is of class $C^2$ on a neighborhood of $0$, then $\Phi_i(f)$ is absolutely continuous with respect to the Lebesgue measure on a neighborhood of $0$ with continuous density. In particular, the density is bounded on some neighborhood of $0$. Since $x\mapsto \zeta(|x|)$ is integrable with respect to the Lebesgue measure and continuous on $\R^n\setminus\{0\}$, the function is integrable with respect to $\Phi_i(f)$.
\end{proof}

Corollary \ref{corollary:RepresentationFunctIntrinsicVolume} naturally leads to the question, under which conditions on $\zeta\in D^n_i$ the integral representation holds for all functions $f\in\Conv(\R^n,\R)$. We will need the following Lemma.
	\begin{lemma}[\cite{ColesantiEtAlHadwigertheoremconvex2020} Lemma 2.15]
	\label{lemma:recover_density}
	Let $1\le i\le n$ and $\zeta\in D^n_i$. For every $t\ge 0$,
	\begin{align*}
		\int_{\R^n}\zeta(|x|)d\Phi_i(u_t,x)=\omega_n\binom{n}{i}\left(t^{n-i}\zeta(t)+(n-i)\int_t^\infty \zeta(s)s^{n-i-1}ds\right).
	\end{align*}
\end{lemma}

\begin{corollary}
	\label{corollary:ZetaIntegrable}
	Let $1\le i\le n-1$ and $\zeta\in D^n_i$. The following are equivalent:
	\begin{enumerate}
		\item $x\mapsto \zeta(|x|)$ is integrable with respect to $\Phi_i(f)$ for all $f\in\Conv(\R^n,\R)$.
		\item $x\mapsto \zeta(|x|)$ is integrable with respect to $\Phi_i(|\cdot|)$
		\item $\int_0^\infty |\zeta(t)|t^{n-i-1}dt<\infty$.
		\item $|\zeta|\in D^n_i$.
	\end{enumerate}
	If any of these condition holds, then $V^*_i(\zeta)[f]=\int_{\R^n}\zeta(|x|)d\Phi_i(f,x)$ for all $f\in\Conv(\R^n,\R)$.
\end{corollary}
\begin{proof}
	$1\Rightarrow 2$ and $3\Rightarrow 4$ are trivial. In order to show $2\Rightarrow 3$, note that we may approximate $|\zeta|$ from below by an increasing sequence of continuous functions, and so dominated convergence and Proposition \ref{lemma:recover_density} imply
	\begin{align*}
		\int_{\R^n}|\zeta(|x|)|d\Phi_i(|\cdot|,x)=\omega_n\binom{n}{i}(n-i)\int_0^\infty |\zeta(s)|s^{n-i-1}ds.
	\end{align*}
	Since the left hand side of this equation is finite, $\int_0^\infty |\zeta(t)|t^{n-i-1}dt<\infty$.\\
	Let us show $4\Rightarrow 1$. Monotone convergence and Theorem \ref{theorem:PrincipalValueRefined} imply that
	\begin{align*}
		\int_{\R^n}|\zeta(x)|d\Phi_i(f,x)=\lim\limits_{\epsilon\rightarrow0}\int_{\R^n\setminus B_\epsilon(0)}|\zeta(x)|d\Phi_i(f,x)=V^*_i(|\zeta|)[f],
	\end{align*}
	where the right hand side is finite since $|\zeta|\in D^n_i$. Thus $x\mapsto \zeta(x)$ is integrable with respect to $\Phi_i(f)$.\\
	The last claim follows directly from Corollary 
	\ref{corollary:RepresentationFunctIntrinsicVolume}.
\end{proof}

	\section{Proof of the Hadwiger theorem on convex functions}
	\label{section:Hadwiger}
	\subsection{Classification in the smooth case}
	\label{section:smoothHadwiger}
	The proof of Theorem \ref{theorem:smoothHadwiger} relies on the classification of all primitive differential $n$-forms on $\R^n\times\R^n$ that are invariant with respect to translations in the second factor. This result requires the following well known facts, which follow easily from L'Hopital's Theorem.
	\begin{lemma}
		\label{lemma:simple_lemma_for_even_differentiable functions}
		If $f:\R\rightarrow\R$ is a smooth even function then there exists $\phi\in C^\infty([0,\infty))$ such that $f(r)=\phi(r^2)$ for all $r\in\R$.
	\end{lemma}
	\begin{lemma}
		\label{lemma:simple_lemma_vanishing_derivatives}
		If $f\in C^\infty(\R)$ satisfies $f^{(i)}(0)=0$ for all $0\le i\le k$ for some $k\in\mathbb{N}$, then there exists $\phi\in C^\infty(\R)$ such that
		\begin{align*}
			f(r)=r^{k+1}\phi(r)\quad\text{for all }r\in\R.
		\end{align*}
	\end{lemma}
	
	\begin{theorem}
		\label{theorem:invariantDiffForm}
		Let $1\le i\le n$. If $\omega\in \Omega^{n-i}(\R^n)\otimes\Lambda^{i}\R^n$ is a primitive $\SO(n)$-invariant differential form, then there exists $c\in\R$ and functions $\phi_1,\phi_2\in C_c^\infty([0,\infty))$ such that
		\begin{align*}
			\omega=c\kappa_i+\phi_1(r^2)\gamma\wedge\tau_i+\phi_2(r^2)\beta\wedge\tau_{i-1}.
		\end{align*}
	\end{theorem}
	\begin{proof}
		Let $\omega\in \Omega^{n-i}(\R^n)\otimes\Lambda^{i}\R^n$ be primitive and $\SO(n)$-invariant. $\omega|_{0}$ defines an $\SO(n)$-invariant differential form in $\Lambda^n(\R^n\times\R^n)$. As the space of primitive elements of this space is generated by the forms $\kappa_l$, $0\le l\le n$, $\omega|_0=c\kappa_i$ for degree reasons. Let us replace $\omega$ by $\omega-c\kappa_i$. We obtain an $\SO(n)$-invariant differential form on $\R^n\times\R^n$ that vanishes in $x=0$. 
		Let us evaluate this form in $x\ne 0$. The stabilizer of this point can be identified with $\mathrm{SO}(n-1)$, and the tangent space splits into the equivariant direct sum
		\begin{align*}
			T_{(x,0)}(\R^n\times\R^n)\cong \R^n\oplus\R^n= (\R x\oplus V)\oplus (\R x\oplus V)\cong (\R x\oplus \R x)\oplus (V\oplus V),
		\end{align*}
		where $V=x^{\perp}$. The stabilizer $\mathrm{Stab}(x)\subset \mathrm{SO}(n)$ of $(x,0)$ acts trivially on $\R x\oplus \R x$ and by the usual diagonal action on $V\oplus V\cong \R^{n-1}\oplus\R^{n-1}$. The evaluation of $\omega$ in $(x,0)$ thus leads to an element of 
		\begin{align*}
			(\Lambda^nT_{(x,0)}(\R^n\times\R^n))^{\mathrm{Stab}(x)}\cong& \bigoplus_{a+b+c=n} \Lambda^a(\R x)\otimes \Lambda^b(\R x)\otimes \Lambda^{c}(\R^{n-1}\oplus \R^{n-1})^{\mathrm{SO}(n-1)},
		\end{align*}
		where $a,b\in\{0,1\}$, $c\in\{n-2,n-1,n\}$. The forms $\gamma$ and $\beta$ span the first two factors, while the last space is generated by the forms $\tau_j$, $0\le j\le n-1$, and the restriction of the symplectic form to $V\oplus V$, where we have used Proposition \ref{proposition:rotation-invariant-forms:SO(n)-invariant-forms_linearized_version} again. Note that the symplectic form of $V\oplus V$ is given by a multiple of $r^2\omega_s-\gamma\wedge\beta$. Thus the space on the right is contained in the space generated by $\gamma,\beta,\omega_s$, and $\tau_j$ for $0\le j\le n-1$. The only $n$-forms are thus the forms $\gamma\wedge\tau_j$, $\beta\wedge\tau_{j-1}$, for $1\le j\le n-1$, as well as multiples of $\omega_s$. The only primitive forms are thus $\gamma\wedge\tau_j$ and $\beta\wedge\tau_{j-1}$, so $\omega$ is a linear combination of these forms. Comparing the degree of homogeneity, we obtain 
		\begin{align*}
			\omega=&\tilde{\phi}_1(x)\gamma\wedge\tau_{i}+\tilde{\phi}_2(x)\beta\wedge\tau_{i-1}
		\end{align*} 
		for two smooth functions $\tilde{\phi}_1,\tilde{\phi}_2:\R^n\setminus\{0\}\rightarrow\R$. Note that $\tilde{\phi}_1$ and $\tilde{\phi}_2$ are uniquely determined by this equation, as the forms $\gamma\wedge\tau_{i}$ and $\beta\wedge\tau_{i-1}$ are linearly independent at each point $x\ne 0$. As $\omega$ is $\SO(n)$-invariant, these functions are thus rotation invariant, that is, we can assume that
		\begin{align*}
			\omega=&\phi_1(r^2)\gamma\wedge\tau_i+\phi_2(r^2)\beta\wedge\tau_{i-1}
		\end{align*} 
		for $r\ne 0$, where $\phi_1,\phi_2\in C^\infty((0,\infty))$. Evaluating these forms along the line $\R(e_1,0)$, we obtain
		\begin{align*}
			&\omega|_{(te_1,0)}\\
			=&\phi_1(t^2)t^2\frac{1}{i!(n-i-1)!}\sum_{\sigma\in S_{n},\sigma(1)=1}\sign(\sigma)dx_{\sigma(1)}\dots dx_{\sigma(n-i)}\wedge dy_{\sigma(n-i+1)}\dots dy_{\sigma(n)}\\
			&+ \phi_2(t^2)t^2 \frac{1}{(i-1)!(n-i)!}\sum_{\substack{\sigma\in S_{n},\\\sigma(n-i+1)=1}}\sign(\sigma)dx_{\sigma(1)}\dots dx_{\sigma(n-i+1)}\wedge  dy_{\sigma(n-i+2)}\dots dy_{\sigma(n)}.
		\end{align*}
		
		The two sums on the right define linearly independent forms on the whole space, so we see that the functions $t\mapsto \phi_j(t^2)t^2$ must be restrictions of smooth functions. As $\omega$ vanishes for $t=0$, the functions
		\begin{align*}
			f_j(r):=\begin{cases}
				\phi_j(r^2)r^2 & \text{for }r\ne 0\\
				0 & \text{for }r=0
			\end{cases}
		\end{align*}
		are thus smooth on $\R$, so we can apply Lemma \ref{lemma:simple_lemma_for_even_differentiable functions} and Lemma \ref{lemma:simple_lemma_vanishing_derivatives} to extend $\phi_1$ and $\phi_2$ to elements of $C^\infty([0,\infty))$. Thus,
		\begin{align*}
			\omega=\phi_1(r^2)\gamma\wedge\tau_i+\phi_2(r^2)\beta\wedge\tau_{i-1}\quad\text{on }\R^n\times\R^n,
		\end{align*}
		which finishes the proof.
	\end{proof}
	The proof of Theorem \ref{theorem:smoothHadwiger} is now a consequence of the following simple relation.
	\begin{lemma}
		\label{lemma:relationDkappa-Betatau}
		For $\phi\in C^\infty_c([0,\infty))$, $0\le i\le n-1$,
		\begin{align*}
			(n-i)\D(\phi(r^2)\kappa_i)=-2\D(\phi'(r^2)\gamma\wedge \tau_i)
		\end{align*}
	\end{lemma}
	\begin{proof}
		Note that $d(r^2)=2\gamma$, so
		\begin{align*}
			d(\phi'(r^2)\gamma\wedge\tau_i)=&-\phi'(r^2)\gamma\wedge d\tau_i=-(n-i)\phi'(r^2)\gamma\wedge\kappa_{i},\\
			d(\phi(r^2)\kappa_{i})=&2\phi'(r^2)\gamma\wedge\kappa_{i}.
		\end{align*}
		As $\D$ vanishes on closed forms, we thus obtain
		\begin{align*}
			0=\D\left(2\phi'(r^2)\gamma\wedge \tau_i+(n-i)\phi(r^2)\kappa_i\right),
		\end{align*}
		which implies the result.
	\end{proof}
	We are now able to prove Theorem \ref{theorem:smoothHadwiger}. Let us restate it in the following form.
	\begin{theorem}
		Let $1\le i\le n$. For every smooth valuation $\mu\in\VConv_i(\R^n)^{\SO(n)}$ there exists $\phi\in C^\infty_c([0,\infty))$ such that
		\begin{align*}
			\mu(f)=D(f)[\phi(|x|^2) \kappa_i]=\binom{n}{i}\int_{\R^n}\phi(|x|^2)d\Phi_i(f,x).
		\end{align*}
	\end{theorem}
	\begin{proof}
		Every smooth $\mu\in\VConv_i(\R^n)^{\SO(n)}$ can be written as $\mu(f)=D(f)[\omega]$ for some $\omega\in \Omega_c^{n-i}(\R^n)\otimes\Lambda^i\R^n$ that is primitive and $\SO(n)$-invariant, compare \cite{KnoerrSmoothvaluationsconvex2024} Corollary 6.7. By Theorem \ref{theorem:invariantDiffForm} there exist $c\in\R$, $\phi_1,\phi_2\in C^\infty([0,\infty))$ such that
		\begin{align*}
			\omega=c\kappa_i+\phi_1(r^2)\gamma\wedge\tau_i+\phi_2(r^2)\beta\wedge\tau_{i-1}.
		\end{align*} 
		As $\omega$ has compact support, we may multiply this expression by $\psi(r^2)$ for some $\psi\in C^\infty_c([0,\infty))$ that is equal to $1$ on $[0,\sqrt{R}]$, where $R>0$ is chosen such that $\supp\omega\subset B_R(0)\times\R^n$. Thus
		\begin{align*}
			\omega=\tilde{\phi}_0(|x|^2)\kappa_i+\tilde{\phi}_1(|x|^2)\gamma\wedge\tau_i+\tilde{\phi}_2(|x|^2)\beta\wedge\tau_{i-1}
		\end{align*}
		for certain $\tilde{\phi}_0,\tilde{\phi}_1,\tilde{\phi}_2\in C^\infty_c([0,\infty))$. 
		Using Proposition \ref{proposition:invariant_forms_Relations}, we obtain
		\begin{align*}
			\D\omega=&\D(\tilde{\phi}_0(|x|^2)\kappa_i)+\D(\tilde{\phi}_1(|x|^2)\gamma\wedge\tau_i)+\D(\tilde{\phi}_2(|x|^2)\beta\wedge\tau_{i-1})\\
			=&\D(\tilde{\phi}_0(|x|^2)\kappa_i)+\D(\tilde{\phi}_1(|x|^2)\gamma\wedge\tau_i)+\D(\tilde{\phi}_2(|x|^2)[|x|^2\kappa_i-\gamma\wedge\tau_{i}).
		\end{align*}
		If we consider
		\begin{align*}
			\Psi_i(t):=-\int_t^\infty \tilde{\phi}_i(s)ds\quad\text{for }i=1,2,
		\end{align*}
		then $\Psi_i\in C_c^\infty([0,\infty))$ and Lemma \ref{lemma:relationDkappa-Betatau} implies
		\begin{align*}
			\D\omega=&\D(\tilde{\phi}_0(|x|^2)\kappa_i)-\frac{n-i}{2}\D(\Psi_1(|x|^2)\kappa_i)+\D(\tilde{\phi}_2(|x|^2)|x|^2\kappa_i)\\
			&+\frac{n-i}{2}\D(\Psi_2(|x|^2)\kappa_i).
		\end{align*}
		Thus the Kernel theorem \ref{theorem:kernel} implies that $\omega$ and 
		\begin{align*}
			\omega':=\left(\tilde{\phi}_0(|x|^2)-\frac{n-i}{2}\Psi_1(|x|^2)+\tilde{\phi}_2(|x|^2)|x|^2+\frac{n-i}{2}\Psi_2(|x|^2)\right)\kappa_i
		\end{align*}
		induce the same valuation, that is,
		\begin{align*}
			\mu(f)=D(f)[\omega']=\binom{n}{i}\int_{\R^n}\phi(|x|^2)d\Phi_i(f,x)
		\end{align*}
		for $\phi\in C^\infty_c([0,\infty))$ given for $t\ge 0$ by
		\begin{align*}
			\phi(t):=\tilde{\phi}_0(t)-\frac{n-i}{2}\Psi_1(t)+\tilde{\phi}_2(t)t+\frac{n-i}{2}\Psi_2(t).
		\end{align*}
	\end{proof}
	Let us remark that one can show that the function $\phi$ is unique by explicitly calculating the symplectic Rumin differential of the form $\phi(|x|^2)\kappa_i$. The uniqueness also follows from the results of the next section, so we refer to \cite{KnoerrSmoothvaluationsconvex2020a} for this calculation.

	\subsection{Proof of Theorem \ref{maintheorem:TopIsomorphism}}
	Our proof of the Hadwiger theorem on convex functions may be seen as a version of the template method, that is, we may recover the "density" $\zeta\in D^n_i$ of the singular Hessian valuation $V^*_{i,\zeta}$ by evaluating this functional in a suitable family of convex functions. The family we will consider is $u_t(x)=\max(0,|x|-t)$, which was also used in the first proof of the Hadwiger theorem on convex functions by Colesanti, Ludwig and Mussnig in \cite{ColesantiEtAlHadwigertheoremconvex2020}.
	The following result is essentially contained in \cite{ColesantiEtAlHadwigertheoremconvex} Lemma 3.7.
	\begin{lemma}
		\label{lemma:inversion_formula}
		Let $1\le i\le n-1$. For $\mu\in\VConv_i(\R^n)$, the function $\zeta:(0,\infty)\rightarrow\R$,
		\begin{align*}
			\zeta(t):=\frac{1}{\omega_n\binom{n}{i}}\left(\frac{\mu(u_t)}{t^{n-i}}+(n-i)\int_t^\infty\frac{\mu(u_s)}{s^{n-i+1}}ds\right),
		\end{align*}
		belongs to $D^n_i$ and satisfies
		\begin{align*}
			\mu(u_t)=\omega_n\binom{n}{i}\left(t^{n-i}\zeta(t)+(n-i)\int_t^\infty \zeta(s)s^{n-i-1}ds\right).
		\end{align*}
		Moreover, $\zeta$ is uniquely determined by the second equation.
	\end{lemma}
	\begin{proof}
		If $R>0$ is such that $\supp\mu\subset B_R(0)$, then $u_t$ coincides with the zero function on a neighborhood of $\supp\mu$ for all $t>R$, so $\mu(u_t)=\mu(0)=0$ for all $t>R$. Thus $t\mapsto \mu(u_t)$ defines a function in $C_c([0,\infty))$, as $t\mapsto u_t$ is a continuous function from $[0,\infty)$ to $\Conv(\R^n,\R)$. It now follows from \cite{ColesantiEtAlHadwigertheoremconvex} Lemma 3.7 that $\zeta\in D^n_i$ and that $\zeta$ is uniquely determined by the second equation.
	\end{proof}
	Our proof of the Hadwiger theorem on convex functions is based on the following observation.
	\begin{lemma}
		\label{lemma:invariant_val_vanishes_ut}
		If $\mu\in\VConv_i(\R^n)^{\SO(n)}$ satisfies $\mu(u_t)=0$ for all $t\ge 0$, then $\mu=0$.
	\end{lemma}
	\begin{proof}
		This is trivial for $i=0$, so let us assume that $1\le i\le n$. Assume that $\mu\in\VConv_i(\R^n)^{\SO(n)}$ satisfies $\mu(u_t)=0$ for all $t\ge0$. There exists a sequence $(\mu_j)_j$ of smooth rotation invariant valuations converging to $\mu$ by Theorem \ref{theorem:approximation_smooth_valuations}. By Proposition \ref{proposition:restriction_support_convergent_sequences}, we may further assume that $\supp\mu,\supp\mu_j\subset B_R(0)$ for all $j\in\mathbb{N}$ for some $R>0$. By Theorem \ref{theorem:approximation_smooth_valuations}, $\mu_j(f)=\int_{\R^n}\zeta_j(|x|)d\Phi_i(f,x)$ for some $\zeta_j\in C_c([0,\infty))$. Moreover, Lemma \ref{lemma:recover_density} shows that $\zeta_j$ is supported on $[0,R]$ for all $j\in\mathbb{N}$, as
		\begin{align}
			\label{eq:inversionFormula}
			\zeta_j(t)=\frac{1}{\omega_n\binom{n}{i}}\left(\frac{\mu_j(u_t)}{t^{n-i}}+(n-i)\int_t^\infty\frac{\mu_j(u_s)}{s^{n-i+1}}ds\right)
		\end{align}
		for all $t>0$. As the map $t\mapsto u_t$ is a continuous map into $\Conv(\R^n,\R)$, the set $\{u_t:t\in [a,b]\}$, $0\le a\le b<\infty$, is compact. In particular, $(\mu_j)_j$ converges uniformly on these sets to $\mu$. For $t>R$, $u_t$ coincides with the zero function on a neighborhood of $B_R(0)$, so $\mu(u_t)=\mu_j(u_t)=0$ for $t>R$. Thus the functions $t\mapsto \mu_j(u_t)$ converge uniformly on $[0,\infty)$ to $t\mapsto \mu(u_t)=0$. Due to \eqref{eq:inversionFormula}, $\zeta_j(t)$ converges uniformly on $[\epsilon,\infty)$ to $0$ for every $\epsilon>0$. We claim that this implies $\supp\mu\subset\{0\}$. \\
		It is enough to show that $\GW(\mu)[\phi_1\otimes\dots\otimes\phi_i]=0$ for all $\phi_1,\dots,\phi_i\in C^\infty_c(\R^n)$ such that $0\notin\supp\phi_k$ for all $1\le k\le i$. As the supports of these functions are compact, there exists $\epsilon>0$ such that $\supp\phi_k\cap B_\epsilon(0)=\emptyset$ for all $1\le k\le i$. It is then easy to see that there exists a constant $C_n$ depending on $n$ only such that
		\begin{align*}
			|\GW(\mu_j)[\phi_1\otimes\dots\otimes\phi_i]|\le C_n\sup_{t\ge\epsilon}|\zeta_j(t)|\cdot \prod_{k=1}^{i}\|\phi_k\|_{C^2(\R^n)}\quad \text{for all }j\in\mathbb{N},
		\end{align*}
		so we can apply Corollary \ref{corollary:weak_continuity_GW} to obtain
		\begin{align*}
			|\GW(\mu)[\phi_1\otimes\dots\otimes\phi_i]|=&\lim\limits_{j\rightarrow\infty}|\GW(\mu_j)[\phi_1\otimes\dots\otimes\phi_i]|\\
			\le& \liminf_{j\rightarrow\infty}C_n \sup_{t\ge\epsilon}|\zeta_j(t)|\cdot \prod_{k=1}^{i}\|\phi_k\|_{C^2(\R^n)}=0,
		\end{align*}
		as $(\zeta_j)_j$ converges uniformly to $0$ on $[\epsilon,\infty)$. Thus $\GW(\mu)[\phi_1\otimes\dots\otimes\phi_i]=0$ unless $0\in\supp\phi_i$ for all $1\le k\le i$. In particular, $\supp\mu\subset\{0\}$, which implies $\mu=0$ due to Proposition \ref{proposition:shape_support}.
	\end{proof}
	Let us remark that Lemma \ref{lemma:invariant_val_vanishes_ut} is equivalent to \cite[Lemma 8.3]{ColesantiEtAlHadwigertheoremconvex2022}, which was established by Colesanti, Ludwig, and Mussnig using Theorem \ref{maintheorem:continuousHadwiger}. As pointed out in the same article, an independent proof of this Lemma is the key step to obtain an alternative proof of Theorem \ref{maintheorem:continuousHadwiger}, which we present now in the proof of Theorem \ref{maintheorem:TopIsomorphism}.
	\begin{proof}[Proof of Theorem \ref{maintheorem:TopIsomorphism}]
		By Corollary \ref{cor:Extension_Vk}, $V^*_i:D^n_{i,R}\rightarrow\VConv_i(\R^n)^{\SO(n)}\cap \VConv_{B_R(0)}(\R^n)$ is continuous and thus uniquely determined by its restriction to the dense subspace $C_c([0,\infty))\cap D^n_{i,R}$. It remains to see that this map is a topological isomorphism.\\
		
		We will construct the inverse map. Let $\mu\in\VConv_i(\R^n)^{\SO(n)}$ and assume that $\supp\mu\subset B_R(0)$. Define $Z(\mu):(0,\infty)\rightarrow\R$ by
		\begin{align*}
			[Z(\mu)](t):=\frac{1}{\omega_n\binom{n}{i}}\left(\frac{\mu(u_t)}{t^{n-i}}+(n-i)\int_t^\infty\frac{\mu(u_s)}{s^{n-i+1}}ds\right).
		\end{align*}
		Then $Z(\mu)\in D^n_i$ by Lemma \ref{lemma:inversion_formula}. Moreover, $u_t$ coincides with the zero function on a neighborhood of $B_R(0)$ for $t>R$. Thus $\mu(u_t)=0$ for $t>R$, which shows $Z(\mu)\in D^n_{i,R}$. We thus obtain a well defined map $Z:\VConv_i(\R^n)^{\SO(n)}\cap\VConv_{B_R(0)}\rightarrow D^n_{i,R}$. By Lemma \ref{lemma:inversion_formula}, for $t\ge 0$
		\begin{align}
			\label{eq:proofMTFormulaDensity}
			\mu(u_t)=\omega_n\binom{n}{i}\left(t^{n-i}[Z(\mu)](t)+(n-i)\int_t^\infty [Z(\mu)](s)s^{n-i-1}ds\right)=V^*_{i,Z(\mu)}(u_t).
		\end{align}
		If $\mu=V^*_i(\zeta)$ for some $\zeta\in D^n_{i,R}$, then  $Z(\mu)=\zeta$ by Lemma \ref{lemma:inversion_formula}, which shows that $Z\circ V^*_i=Id_{D^n_{i,R}}$. In particular, $V^*_i$ is injective.\\
		Let us show that it is also surjective. As $\mu(u_t)=V^*_{i,Z(\mu)}(u_t)$ for all $t\ge 0$, the valuation $\tilde{\mu}:=\mu-V^*_{i,Z(\mu)}\in\VConv_i(\R^n)^{\SO(n)}$ satisfies $\tilde{\mu}(u_t)=0$ for all $t\ge0$. Lemma \ref{lemma:invariant_val_vanishes_ut} implies $\tilde{\mu}=0$, that is, $\mu=V^*_{i,Z(\mu)}=V^*_i(Z(\mu))$. Thus $V^*_i$ is also surjective and $Z=(V^*_i)^{-1}$\\
		To see that $(V^*_i)^{-1}$ is also continuous, note that $D^n_{i,R}$ and $\VConv_i(\R^n)^{\SO(n)}\cap \VConv_{B_R(0)}(\R^n)$ are both Banach spaces. Thus $(V^*_i)^{-1}$ is continuous by the open mapping theorem.
	\end{proof}
	
	\paragraph*{Acknowledgments}
	I want to thank Fabian Mussnig for his comments on the first draft of this article, which lead to the representation of the functional intrinsic volumes as principal value integrals.

	\bibliography{../../../library/library.bib}
	\Addresses

\end{document}